\documentclass[11pt]{amsart}
\usepackage{amssymb,amsmath,amsfonts,color,hyperref}
\usepackage[all,cmtip]{xy}

\newtheorem{proposition}{Proposition}[section]
\newtheorem{definition}[proposition]{Definition}
\newtheorem{lemma}[proposition]{Lemma}
\newtheorem{theorem}[proposition]{Theorem}
\newtheorem{corollary}[proposition]{Corollary}

\newtheorem{maintheorem}{Theorem}
\theoremstyle{definition}
\newtheorem*{remark}{Remark}

\DeclareMathOperator{\vol}{vol}

\DeclareMathOperator{\Gr}{Gr}

\DeclareMathOperator{\SO}{SO}

\DeclareMathOperator{\nc}{nc}
\DeclareMathOperator{\n}{n}

\DeclareMathOperator{\FlagArea}{FlagArea}

\DeclareMathOperator{\inte}{\mathrm{int}}

\renewcommand{\O}{\mathrm{O}}
\newcommand{\spann}{\mathrm{span}}

\newcommand{\R}{\mathbb{R}}

\newcommand{\N}{\mathbb{N}}

\newcommand{\K}{\mathcal{K}}

\newcommand{\Id}{\mathrm{Id}}

\newcommand{\F}{\mathcal{F}}

\newcommand{\W}{W}
\newcommand{\V}{V}

\newcommand\flag[2]{\left[\begin{array}{c} #1\\ #2 \end{array}
  \right]} 
\newcommand{\largewedge}{\mbox{\Large $\wedge$}}

\newtheorem{note}{Remark} 
\def\note#1{\ifvmode\leavevmode\fi\vadjust{\vbox to0pt{\vss
 \hbox to 0pt{\hskip\hsize\hskip1em
\vbox{\hsize3.5cm\small\raggedright\pretolerance10000
 \noindent #1\hfill}\hss}\vbox to8pt{\vfil}\vss}}}
\renewcommand{\:}{\, : \,}

\title{Flag area measures}

\author{Judit Abardia-Ev\'equoz} 
\address{Institut f\"ur Mathematik, Goethe-Universit\"at Frankfurt am Main, 
Robert-Mayer-Str. 10, 60054 Frankfurt, Germany}
\email{abardia@math.uni-frankfurt.de}

\author{Andreas Bernig} 
\address{Institut f\"ur Mathematik, Goethe-Universit\"at Frankfurt am Main, 
Robert-Mayer-Str. 10, 60054 Frankfurt, Germany}
\email{bernig@math.uni-frankfurt.de}

\author{Susanna Dann} 
\address{Departamento de Matem\'aticas, Universidad de los Andes, Carrera 1 N\textsuperscript{\underline{o}} 18A 12, 111711 Bogot{\'a}, Colombia}
\email{s.dann@uniandes.edu.co}

\begin{document}



\begin{abstract}
A flag area measure on an $n$-dimensional euclidean vector space is a continuous translation-invariant valuation with values in the space of signed measures on the flag manifold consisting of a unit vector $v$ and a $(p+1)$-dimensional linear subspace containing $v$ with $0 \leq p \leq n-1$. 

Using local parallel sets, Hinderer constructed examples of $\mathrm{SO}(n)$-covariant flag area measures. There is an explicit formula for his flag area measures evaluated on polytopes, which involves the squared cosine of the angle between two subspaces.

We construct a more general sequence of smooth $\mathrm{SO}(n)$-covariant flag area measures via integration over the normal cycle of appropriate differential forms. We provide an explicit description of our measures on polytopes, which involves an arbitrary elementary symmetric polynomial in the squared cosines of the principal angles between two subspaces. 

Moreover, we show that these flag area measures span the space of all smooth $\mathrm{SO}(n)$-covariant flag area measures, which gives a classification result in the spirit of Hadwiger's theorem. 
\end{abstract}

\maketitle
\section{Introduction}

\subsection{General background}

Let $(\V,\langle\cdot,\cdot\rangle)$ be a euclidean vector space of dimension $n$ and let $A$ be an abelian semigroup. Let $\mathcal{K}(\V)$ be the space of compact convex subsets of $\V$. A {\it valuation} on $\V$ is a map $\mu:\mathcal{K}(\V) \to A$ satisfying  
\begin{displaymath}
 \mu(K \cup L)+\mu(K \cap L)=\mu(K)+\mu(L),
\end{displaymath}
whenever $K,L,K \cup L \in \mathcal{K}(\V)$. If $A$ is a topological abelian semigroup, then $\mu$ is a {\it continuous} valuation if it is continuous with respect to the topology on $\mathcal{K}(\V)$ induced by the Hausdorff metric.

Natural properties of valuations are translation-invariance (i.e. $\mu(K+v)=\mu(K)$ for all $K$ and $v \in \V$) and rotation-invariance (i.e. $\mu(g K)=\mu(K)$ for all $K$ and $g \in \mathrm{SO}(n)$). Hadwiger's seminal theorem states that the space of real-valued, continuous, rotation- and translation-invariant valuations is spanned by the {\it intrinsic volumes} $\mu_0,\ldots,\mu_n$. We refer to \cite{klain_rota} for a proof of this theorem and a detailed account on the related kinematic formulas in euclidean vector spaces. 

Intrinsic volumes may be defined via a tube formula (called Steiner's formula, see \cite{steiner, klain_rota,schneider_book14}) as follows. Denote by $B^m_2$ the euclidean unit ball of dimension $m$ and by $\kappa_{m}:=\frac{\pi^{\frac{m}{2}}}{\Gamma(1+\frac{m}{2})}$ its volume.  Then $\omega_m:=m\kappa_m$ is the volume of the unit sphere of dimension $m-1$, denoted by $S^{m-1}$. Given $K \in \mathcal{K}(\V)$ and $\rho >0$, let $K+\rho B^n_2$ be the Minkowski sum of $K$ and a $\rho$-ball, which we call a {\it $\rho$-tube around $K$}. Then $\vol(K+\rho B^n_2)$ is a polynomial in $\rho$, whose coefficients (up to a normalization) are the intrinsic volumes:
\begin{displaymath}
 \vol(K+\rho B^n_2)=\sum_{k=0}^n \mu_k(K) \kappa_{n-k} \rho^{n-k}. 
\end{displaymath}

Intrinsic volumes appear in the following kinematic formulas due to Chern-Blaschke-Santal\'o. Let $\overline{\mathrm{SO}(n)} := \V \rtimes \mathrm{SO}(n)$ be the special euclidean group (i.e. the group of isometries of $V$ preserving orientation), endowed with the product of the Lebesgue measure on $\V$ and the probability measure on ${\mathrm{SO}(n)}$. Given $K,L\in \mathcal{K}(\V)$, the {\it intersectional kinematic formulas} assert that for $0\leq i \leq n$ (see e.g. \cite{klain_rota,schneider_book14})
\begin{equation} \label{eq_intersectional_kinform}
 \int_{\overline{\mathrm{SO}(n)}} \mu_i(K \cap \bar g L) d\bar g = \flag{n+i}{i} \sum_{k+l=n+i}  \flag{n+i}{k}^{-1} \mu_k(K) \mu_l(L).
\end{equation}
The expression on the right-hand side involves {\it flag coefficients} defined by 
\begin{displaymath}
 \flag{n}{k}:=\binom{n}{k} \frac{\kappa_n}{\kappa_k \kappa_{n-k}}.
\end{displaymath}
The {\it additive kinematic formulas}, also called {\it rotation sum} formulas (see \cite{schneider_book14}), are given by
\begin{equation} \label{eq_additive_kinform}
 \int_{\mathrm{SO}(n)} \mu_i(K+gL)dg=\flag{2n-i}{n-i} \sum_{k+l=i}  \flag{2n-i}{n-k}^{-1} \mu_k(K) \mu_l(L),
\end{equation}
with $0\leq i \leq n$. 

Both formulas admit local versions. In order to describe them, let us recall the notion of a support measure (we refer to \cite{schneider_book14} for a detailed study). Let $SV:=\V \times S^{n-1}$ be the unit sphere bundle over $\V$. Denote by $\mathcal{B}(SV)$ the Borel subsets of $SV$. For a given $\eta \in \mathcal{B}(SV)$, consider the part of the $\rho$-tube around $K$ defined by 
\begin{displaymath}
M_\rho(K,\eta):=\{x \in \V \: 0<d(x,K)\leq \rho; (p(K,x),u(K,x)) \in \eta\}.
\end{displaymath}
Here $p(K,x)$ is the unique nearest point to $x$ in $K$ and $u(K,x):=\frac{x-p(K,x)}{\|x-p(K,x)\|}$ is a unit normal vector to $K$ at $p(K,x)$. The volume of $M_\rho(K,\eta)$ is again a polynomial in $\rho$, whose coefficients (up to a normalization) are the {\it support measures} of $K$:
\begin{displaymath}
 \vol M_\rho(K,\eta)=\frac{1}{n} \sum_{k=0}^{n-1} \Theta_k(K,\eta) \binom{n}{k} \rho^{n-k}.  
\end{displaymath}
For a fixed $K$, $\Theta_k(K,\cdot)$ is a measure on $SV$. For a fixed $\eta$, $\Theta_k(\cdot,\eta)$ is a measurable valuation. Moreover, $K \mapsto \Theta_k(K,\cdot)$ is continuous, where the space of measures is endowed with the weak topology. 

The support measures $\Theta_k$ are $\SO(n)$-invariant, i.e. $\Theta_k(gK,g\eta)=\Theta_k(K,\eta)$ for every $g\in \SO(n), K\in\K(V)$ and $\eta\in\mathcal{B}(SV)$. For characterization theorems in the spirit of Hadwiger, see \cite{schneider78} and \cite[p. 226]{schneider_book14}. The marginals of support measures are {\it curvature measures} and {\it area measures}:
\begin{align*}
 C_k(K,\beta) & := \Theta_k(K, \beta \times S^{n-1}), \quad \beta \in \mathcal{B}(\V),\\
 S_k(K,\beta) & := \Theta_k(K, \V \times \beta), \quad \beta \in \mathcal{B}(S^{n-1}), 
\end{align*}
where $0 \leq k \leq n-1$. Another common normalization is the following: 
\begin{align*}
 \Phi_k(K,\cdot) & := \frac{1}{n \kappa_{n-k}}\binom{n}{k} C_k(K,\cdot),\\
 \Psi_k(K,\cdot) & := \frac{1}{n \kappa_{n-k}}\binom{n}{k} S_k(K,\cdot).
\end{align*}
One completes the definition by setting $\Phi_n(K,\beta):=\mathcal{H}^n(K \cap \beta)$. Then $\Phi_k(K,\V)=\mu_k(K)$ for $0 \leq k \leq n$ and $\Psi_k(K,S^{n-1})=\mu_k(K)$ for $0 \leq k \leq n-1$.

The local form of the formula \eqref{eq_intersectional_kinform} can best be described in terms of the curvature measures $\Phi_k$. For $\beta,\beta' \in \mathcal{B}(\V)$, we have
\begin{equation} \label{eq_local_intersectional_kinform}
 \int_{\overline{\mathrm{SO}(n)}} \Phi_i(K \cap \bar g L, \beta \cap \bar g \beta') d\bar g = \flag{n+i}{i} \sum_{k+l=n+i}  \flag{n+i}{k}^{-1} \Phi_k(K,\beta) \Phi_l(L,\beta'),
\end{equation}
with $0\leq i\leq n$.
The additive kinematic formula \eqref{eq_additive_kinform} admits a localization in terms of the area measures $S_k$ as follows, for $\beta, \beta' \in \mathcal{B}(S^{n-1})$
\begin{displaymath}
\int_{\mathrm{SO}(n)} S_i(K+gL,\beta \cap g\beta') dg=\frac{1}{\omega_n} \sum_{k+l=i} \binom{i}{k} S_k(K, \beta) S_l(L,\beta'), \quad 0 \leq i \leq n-1. 
\end{displaymath}

Formulas of this type can also be proved for subgroups of the euclidean motion group $\overline{\mathrm{SO}(n)}$. For the important hermitian case, where the rotation group is replaced by the unitary rotation group $\overline{\mathrm{U}(m)}$, the global kinematic formulas were obtained in \cite{bernig_fu_hig}; the local intersectional kinematic formulas were derived in \cite{bernig_fu_solanes}; and the local additive formulas were found in \cite{wannerer_area_measures}. 

Let us next recall the construction of flag-type support, area and curvature measures via local tube formulas. We refer to \cite{hug_tuerk_weil} for a survey and to \cite{ goodey_hinderer_hug_rataj_weil,hinderer_hug_weil, hug_rataj_weil} for more recent developments. Let $\Gr_p(\V)$ denote the Grassmann variety of $p$-dimensional subspaces in $\V$ and $\overline{\Gr}_p(\V)$ the set of all affine $p$-planes in $\V$. Then $\Gr_p(\V)$ admits an $\mathrm{SO}(n)$-invariant measure, while $\overline{\Gr}_p(\V)$ admits an $\overline{\mathrm{SO}(n)}$-invariant measure. These measures are unique up to a normalization, and we choose the standard normalization from \cite[Thm. 13.2.12]{schneider_weil08}. For a given $\bar E \in \overline{\Gr}_p(\V)$, denote by $E \in \Gr_p(\V)$ the linear space parallel to $\bar E$. 

Fix $K \in \mathcal{K}(\V)$ and $\rho>0$. Let $\bar E \in \overline{\Gr}_p(\V)$ be such that $\bar E \cap K = \emptyset$. For almost all such $\bar E$, there exists a unique pair of nearest points $p(K,\bar E) \in K$ and $l(K,\bar E) \in \bar E$. Then $u(K,\bar E):=\frac{l(K,\bar E)-p(K,\bar E)}{d(K,\bar E)} \in S^{n-1}$ is a normal vector to $K$ at $p(K,\bar E)$. Given a Borel subset $\eta$ of $\V \times S^{n-1} \times \Gr_p(\V)$, the local parallel set of $K$ in $\overline{\Gr}_p(\V)$ is defined by 
\begin{displaymath}
 M_\rho^{(p)}(K,\eta):=\{\bar E \in \overline{\Gr}_p(\V) \: 0 < d(K,\bar E) \leq \rho, (p(K,\bar E),u(K,\bar E),E) \in \eta\}. 
\end{displaymath}
The volume of $M_\rho^{(p)}(K,\eta) \subset \overline \Gr_p(V)$, $0\leq p\leq n-1$, is a polynomial in $\rho$, 
\begin{displaymath}
 \vol M_\rho^{(p)}(K,\eta)= \sum_{k=0}^{n-p-1} \Xi^{(p)}_k(K,\eta) \kappa_{n-p-k} \rho^{n-p-k},
\end{displaymath}
whose coefficients $\Xi^{(p)}_k(K,\cdot), 0 \leq k \leq n-p-1$ are the {\it flag-type support measures} (see \cite{hug_tuerk_weil}). Another normalization given in \cite{hug_tuerk_weil} is the following:
\begin{displaymath}
 \Theta_k^{(p)} := \omega_{n-p-k} \binom{n-p-1}{k}^{-1} \Xi_k^{(p)}.
\end{displaymath}
Let us mention some properties of flag-type support measures. For a fixed $\eta$, the map $\Xi^{(p)}_k(\cdot,\eta)$ is a measurable valuation. The map $K \mapsto \Xi^{(p)}_k(K,\cdot)$ is continuous. The group $\overline{\mathrm{SO}(n)}$ acts on $\V \times S^{n-1} \times \Gr_p(\V)$ by $\overline g\cdot(x,v,E)=(gx+t,gv,gE)$, where $\overline g=(g,t)$ with  $g\in\SO(n)$ and $t \in V$. The maps $\Xi^{(p)}_k$ are $\overline{\mathrm{SO}(n)}$-invariant by construction, i.e. for $\bar g \in \overline{\mathrm{SO}(n)}$, $\Xi^{(p)}_k(\bar gK,\bar g\eta)=\Xi^{(p)}_k(K,\eta)$.

Similarly to curvature and area measures, which appear as marginals of support measures, {\it flag-type curvature} and {\it flag-type area measures} are defined by
\begin{align*}
	\Phi_k^{(p)}(K,\beta) & := \Xi^{(p)}_k(K,\beta \times S^{n-1}), \quad \beta \in \mathcal{B}(\V \times \Gr_p(\V)),\\
	\Psi_k^{(p)}(K,\beta) & := \Xi^{(p)}_k(K,\V \times \beta), \quad \beta \in \mathcal{B}(S^{n-1} \times \Gr_p(\V)).
\end{align*}
We will use the following normalization of flag-type area measures, as used in \cite{hinderer_hug_weil}: 
\begin{displaymath}
	S_k^{(p)} := \omega_{n-p-k} \binom{n-p-1}{k}^{-1} \Psi_k^{(p)}.
\end{displaymath}
Set
\begin{align*}
	F(n, p+1) & := \{(v,E) \in S^{n-1} \times \Gr_{p+1}(\V) \: v \in E\}, \\
	F^\perp(n,p) & := \{(v,E) \in S^{n-1} \times \Gr_p(\V) \: v \perp E\}.
\end{align*}
The diffeomorphisms $F^\perp(n,p) \cong F(n,p+1), (v,E) \mapsto (v,\mathbb R v \oplus E)$  and $F(n,p+1) \cong F(n,n-p), (v,E) \mapsto (v,\R v \oplus E^\perp)$ are $\O(n)$-equivariant. By construction the measure $S_k^{(p)}(K, \cdot)$ is concentrated on $F^\perp(n,p)$.

We now give an explicit expression of flag-type area measures on polytopes, which follows from \cite{hinderer_hug_weil}. Given a polytope $P$, denote by $\mathcal{F}_k(P)$ the set of all $k$-dimensional faces of $P$. For a face $F$, let $N(P,F)$ denote the normal cone of $P$ at $F$ and let $\n(P,F):=N(P,F) \cap S^{n-1}$. Let $F$ also denote the linear space parallel to the face $F$ of the same dimension. The cosine of the angle between two subspaces $E$ and $F$ is denoted by $|\cos(E,F)|$, see Section \ref{sec_jordan} for its definition.

\begin{proposition}[Theorem 3.8 in \cite{hinderer_hug_weil}] \label{prop_explicit_smk}
 Let $P$ be a polytope, $0\leq p \leq n-1$, $0\leq k \leq n-p-1$ and $\beta \in \mathcal{B}(F^\perp(n,p))$. Then 
 \begin{align*}
  S_k^{(p)}(P,\beta) 
  & = \binom{n-p-1}{k}^{-1} \frac{\omega_{n-p}}{\omega_n} \times \\
  & \times \sum_{F \in \mathcal{F}_k(P)} \vol_k(F) \int_{\n(P,F)} \int_{\Gr_{p+1}(v)} \mathbf{1}_{(v,E \cap v^\perp) \in \beta} \cos^2(E^\perp,F) dE \ dv.
 \end{align*}
Here $\Gr_{p+1}(v)$ denotes the Grassmannian of all $(p+1)$-planes containing $v$, endowed with an invariant probability measure $dE$.
\end{proposition}

These flag area measures appear naturally in several contexts: extension of valuations \cite{hinderer_hug_weil}; integral formulas for projection functions of convex bodies \cite{goodey_hinderer_hug_rataj_weil}; integral formulas for mixed volumes of convex bodies \cite{hug_rataj_weil}.

\subsection{Results of the present paper}

We start with a general definition.

\begin{definition}
Let $0 \leq p \leq n-1$. A \emph{flag area measure} on $V$ is a translation-invariant valuation $\Phi$ with values in the space of signed measures on the flag manifold $F(n,p+1)$. The space of continuous flag area measures is denoted by $\FlagArea^{(p)}$ and its $k$-homogeneous elements by $\FlagArea^{(p)}_k$. For a group $G$ acting linearly on $V$, $\Phi$ is called \emph{covariant} if $\Phi(gK)=g_*\Phi(K)$ for all $g \in G$. The subspace of $G$-covariant continuous flag area measures is denoted by $\FlagArea^{(p),G}$.
\end{definition}

Without further assumptions, like for example a version of the local definedness from \cite[Satz 2]{schneider75b} or \cite[Theorem 6.1]{schneider78}, this definition is probably too general to be useful. In Section \ref{sec_smooth} we introduce a notion of smoothness of flag area measures, which is stronger than continuity. The main purpose of this article is to classify $\SO(n)$-covariant smooth flag area measures. 

To describe our results more precisely, we need the notion of Jordan angles between subspaces. Let $\W$ be a euclidean vector space and $0 \leq k,p \leq \dim \W$. The orthogonal group $\mathrm{O}(\W)$ acts diagonally on the product $\Gr_p(\W) \times \Gr_k(\W)$. The orbits under this action are characterized in terms of Jordan angles $\theta_1,\ldots,\theta_m$, where $m=\min\{k,\dim \W-k,p,\dim \W-p\}$. Given $E \in \Gr_p(\W)$ and $F \in \Gr_k(\W)$, we denote by
$$\sigma_i(E,F):=\sigma_i(\cos^2\theta_1,\ldots,\cos^2 \theta_m)$$ the $i$-th elementary symmetric polynomial. For instance, $\sigma_m(E,F)$ is the square of the usual cosine between $E$ and $F$. 

In the special case when $\dim \W=2a$ is even and $E,F$ are oriented subspaces with $\dim E=\dim F=a$, there is an additional $\SO(W)$-invariant, $\tilde\sigma_a(E,F)$. 
We refer to Section \ref{sec_jordan} for the definition and  properties of Jordan angles and of $\tilde \sigma_a$.

As our first main result, we construct a $3$-parameter family of flag area measures $S_k^{(p),i}$, which contains the $2$-parameter family $S_k^{(p)}$ as a subfamily.

\begin{maintheorem}[Construction of flag area measures] \label{mainthm_construction}
For every $0 \leq p,k \leq n-1, 0 \leq i \leq m:=\min\{k,n-k-1,p,n-p-1\}$, there exists a unique continuous translation-invariant flag area measure such that for a polytope $P \subset \V$ and $\beta \in \mathcal{B}(F(n,p+1))$,
\begin{align}
  S_k^{(p),i}(P,\beta) & =c_{n,k,p,i} \sum_{F \in \mathcal{F}_k(P)} \vol_k(F) \int_{\n(P,F)} \int_{\Gr_{p+1}(v)} \mathbf{1}_{(v,E) \in \beta} \sigma_i(E^\perp,F) dE \ dv, \label{eq_Skpi}
\end{align}
 where 
\begin{displaymath}
c_{n,k,p,i}:=\binom{n-1}{k}^{-1} \binom{m}{i}^{-1} \binom{|k-(n-1-p)|+m}{i}^{-1} \binom{n-1}{i}.
\end{displaymath}
The Jordan angles are computed with respect to the $(n-1)$-dimensional space $\W:=v^\perp$. 
 
For an odd $n$, there exists an additional unique continuous translation-invariant flag area measure $\tilde S_\frac{n-1}{2}^{(\frac{n-1}{2})}$ such that for a polytope $P \subset \V$ and $\beta \in \mathcal{B}(F(n,\frac{n+1}{2}))$,
\begin{equation}
 \tilde S_\frac{n-1}{2}^{(\frac{n-1}{2})}(P,\beta) = \sum_{F \in \mathcal{F}_\frac{n-1}{2}(P)}\!\! \vol_{\frac{n-1}{2}}(F) \int_{\n(P,F)} \int_{\Gr_{\frac{n+1}{2}}(v)} \!\!\mathbf{1}_{(v,E) \in \beta}\, \tilde\sigma_{\frac{n-1}{2}}(E^{\perp},F) dE \ dv. \label{eq_Sk_sp}
 \end{equation} 
\end{maintheorem}

The flag area measures from Theorem \ref{mainthm_construction} can also be computed explicitly for a smooth compact convex body $K$. Let $x \in \partial K$ be a boundary point and let $E \in \Gr_{p+1}(V)$ be a plane containing the normal vector $\nu(x)$ to $K$ at the point $x$. There are three linear maps on the tangent space $T_x \partial K$: the shape operator $S_x$ and the orthogonal projections $\Pi_E$ and $\Pi_{E^\perp}$ onto $E \cap T_x\partial K$ and $E^\perp$, respectively. We orient the spaces $E \cap T_x\partial K$ and $E^\perp$ in such a way that $\R \nu(x) \oplus (E \cap T_x\partial K) \oplus E^\perp \cong V$ as oriented vector spaces. Let $D$ denote the mixed discriminant. We derive the following formulas.  

\begin{maintheorem} \label{mainthm_smooth}
Let $k,p,i,\beta$ be as in Theorem \ref{mainthm_construction}. Then for a smooth compact convex body $K$,
\begin{align}
& S_k^{(p),i}(K,\beta)  = c_{n,k,p,i} \binom{n-1}{k} \sum_{a=\min\{k,p\}-m}^{\min\{k,p\}-i} \binom{\min\{k,p\}-a}{i} \binom{k}{a} \times \nonumber \\ 
& \times \int_{\partial K} \int_{\Gr_{p+1}(\nu(x))} \mathbf{1}_{(\nu(x),E) \in \beta} D(S_x[n-k-1],\Pi_E[a],\Pi_{E^\perp}[k-a]) dE d\mathcal{H}^{n-1}(x). \label{eq_flagOnSmooth}
\end{align}
\noindent
If $n$ is odd and $p=k=\frac{n-1}{2}$,  
\begin{align*}
 &\tilde S^{\left(\frac{n-1}{2}\right)}_\frac{n-1}{2}(K,\beta) = \\ 
 &(-1)^\frac{n-1}{2} \int_{\partial K} \int_{\Gr_{p+1}(\nu(x))} \mathbf{1}_{(\nu(x),E) \in \beta} \det(\Pi_{E^\perp} \circ S_x:E \cap T_x\partial K \to E^\perp) dE d\mathcal{H}^{n-1}(x).
\end{align*}
\end{maintheorem} 

\begin{remark}
\begin{enumerate}
\item We will see that in the case $k \leq p$, \eqref{eq_flagOnSmooth} simplifies to 
\begin{align*}
& S_k^{(p),i}(K,\beta) = c_{n,k,p,i} \binom{n-1}{k} \binom{k}{i} \times \\ 
& \quad \times \int_{\partial K} \int_{\Gr_{p+1}(\nu(x))} \mathbf{1}_{(\nu(x),E) \in \beta} D(S_x[n-k-1],\mathrm{Id}[k-i],\Pi_{E^\perp}[i])dE d\mathcal{H}^{n-1}(x).
\end{align*}
\item In Section \ref{sec_construction} we will prove an integral formula which generalizes \eqref{eq_flagOnSmooth} to arbitrary convex bodies. It involves generalized principal curvatures and the generalized shape operator. 
\end{enumerate}
\end{remark}
	
The next theorem summarizes the main properties of flag area measures. Properties (i)-(v) are easy to prove. The proof of (vi) requires some deep results by James about the distribution of Jordan angles and by Aomoto about Selberg-type integrals.

\begin{maintheorem}[Properties of flag area measures] \label{mainthm_properties}
The flag area measures $S_k^{(p),i}$ and $\tilde S_\frac{n-1}{2}^{(\frac{n-1}{2})}$ satisfy: 
\begin{enumerate}
	\item[(i)] For $k \leq n-p-1$, we have $m=\min\{k,p\}$ and $S_k^{(p),m}=\frac{\omega_{n}}{\omega_{n-p}}S_k^{(p)}$.
	\item[(ii)] For a fixed $\beta$, $S_k^{(p),i}(\cdot,\beta)$ and $\tilde S_\frac{n-1}{2}^{(\frac{n-1}{2})}(\cdot,\beta)$ are translation-invariant valuations homogeneous of degree $k$ and $\frac{n-1}{2}$ respectively. 
 	\item[(iii)] For all $g \in \mathrm{O}(n)$, $S_k^{(p),i}(gK,g\beta)=S_k^{(p),i}(K,\beta)$.
 	\item[(iv)] For every $K\in \mathcal{K}(V)$, the measures $S_k^{(p),i}(K, \cdot)$ are positive.
 	\item[(v)] For all $g \in \mathrm{O}(n)$, $\tilde S_\frac{n-1}{2}^{(\frac{n-1}{2})}(gK,g\beta)=\det g \cdot \tilde S_\frac{n-1}{2}^{(\frac{n-1}{2})}(K,\beta)$.
 	\item[(vi)] Let $\pi:F(n,p+1) \to S^{n-1}$ be the projection onto the first factor. For each $\beta \in \mathcal{B}(S^{n-1})$, we have 
 \begin{displaymath}
  S_k^{(p),i}(K,\pi^{-1}(\beta))= S_k(K,\beta) \quad \text{ and } \quad \tilde S_\frac{n-1}{2}^{(\frac{n-1}{2})}(K,\pi^{-1}(\beta))=0.
 \end{displaymath}
\end{enumerate}
\end{maintheorem}

To prove Theorem \ref{mainthm_construction} we construct a sequence of $\overline{\mathrm{SO}(n)}$-invariant differential forms on the product $\V \times F(n,p+1)$. The flag area measures are obtained by integration over the pull-back of the normal cycle of $K$ under the projection map $\V \times F(n,p+1) \to \V \times S^{n-1}, (x,v,E) \mapsto (x,v)$. The flag area measures obtained in this way will be called smooth and the corresponding space will be denoted by $\mathrm{FlagArea}^{(p),sm}$ (cf. Definition \ref{def:smoothFlagMeas}).

\begin{maintheorem} \label{thm_basis}
Let $0 \leq p,k \leq n-1$. For $(p,k) \neq \left(\frac{n-1}{2},\frac{n-1}{2}\right)$ a basis of $\FlagArea^{(p),sm,\mathrm{SO}(n)}_k$ is given by:
$$ S_k^{(p),i} \text{ with } 0 \leq i \leq m=\min\{k,n-k-1,p,n-p-1\}.$$
For $n$ odd and $(p,k) = \left(\frac{n-1}{2},\frac{n-1}{2}\right)$, a basis of $\FlagArea^{\left(\frac{n-1}{2}\right),sm,\mathrm{SO}(n)}_{\frac{n-1}{2}}$ 
is given by 
\begin{displaymath}
 S_\frac{n-1}{2}^{\left(\frac{n-1}{2}\right),i} \text{ with } 0 \leq i \leq \frac{n-1}{2}, \quad \text{ and } \quad \tilde S_\frac{n-1}{2}^{\left(\frac{n-1}{2}\right)}.
\end{displaymath}
\end{maintheorem}

Our construction of flag area measures can be generalized in several ways. First, one could use other - partial or complete - flag manifolds instead of $F(n,p+1)$. 
Second, the $\SO(n)$-invariance may be relaxed or dropped. This will lead to a very general class of smooth flag area measures that are translation-invariant, but not necessarily $\SO(n)$-invariant. We think that these general flag area measures deserve further study.

\subsection*{Plan of the paper}

In Section \ref{sec_smooth}, we introduce the notion of smooth flag area measures. These are given by forms on $V\times F(n,p+1)$. We describe the forms on this space that induce the trivial flag area measure. 

In Section \ref{sec_dimension}, we determine the dimension of the space of $\SO(n)$-covariant smooth flag area measures that are homogeneous of degree $k$, $0\leq k\leq n-1$.

In Section \ref{sec_jordan}, we recall the notion of Jordan angles between two subspaces of a euclidean vector space and compute the average value of elementary symmetric functions over the squared cosines of the Jordan angles between two subspaces. 

The main results of this paper are proved in Section \ref{sec_construction}. We give an explicit construction of the differential forms inducing linearly independent flag area measures. We then write the value of our flag area measures on a compact convex body in terms of generalized principal curvatures and the generalized shape operator. Theorems \ref{mainthm_construction} and \ref{mainthm_smooth} are easy consequences of this general formula.

\section{Smooth flag area measures}\label{sec_smooth}

In the following, given a smooth manifold $M$, we denote the space of differential forms on $M$ by $\Omega^*(M)$, the tangent (resp. cotangent) bundle by $TM$ (resp. $T^*M$) and for $p\in M$, the tangent (resp. cotangent) space of $M$ at $p$, by $T_pM$ (resp. $T_p^*M$). 

Let us recall the definition of the {\it fiber integration}, also called the {\it push-forward} of differential forms. We follow the sign convention in  \cite{alvarez_fernandes}. For another sign convention see e.g.\! \cite{berline_getzler_vergne}. 

\begin{definition}\label{def_push_forward}
Let $M$ and $B$ be oriented manifolds. Let $\Pi:M\to B$ be a fiber bundle with a compact fiber of dimension $r$ oriented by the local product orientation. Let $d\geq 0$. The \emph{fiber integration} of a form $\eta\in\Omega^{d+r}(M)$ of degree $d+r$ is the form $\Pi_*\eta \in\Omega^{d}(B)$ of degree $d$ defined by 
\begin{equation}\label{eq_fiber_int_gen}
 \Pi_*\eta|_{y}(w_1,\dots,w_{d}):=\int_{\Pi^{-1}(y)}\eta_{w_1,\dots,w_{d}}, \quad y \in B,
\end{equation}
where $\eta_{w_1,\dots,w_{d}} \in\Omega^{r}(M)$ is defined as follows: for $x\in M$ with $\Pi(x)=y$ and for $W_j\in T_{x}M$ with $d\Pi(W_j)=w_j$,
\begin{displaymath}
(\eta_{w_1,\dots,w_d})|_{x}:=\eta|_{x}(W_1,\dots,W_{d},-). 
\end{displaymath}
$W_j$ are called \emph{lifts of $w_j$}. The lifts $W_j$ are not unique. However, the right-hand side of \eqref{eq_fiber_int_gen} is independent of their choices.
\end{definition}

An equivalent definition is as follows. Let $\Pi^*\beta$ denote the pullback of a form $\beta$ under the map $\Pi$. It has the same degree as $\beta$. The fiber integration, $\Pi_*\eta$, is uniquely defined by the equation
\begin{displaymath}
 \int_B\Pi_*\eta \wedge\beta:=(-1)^{r(\dim B-d)} \int_M\eta\wedge\Pi^*\beta \, ,
\end{displaymath}
to be satisfied for every compactly supported differential form $\beta\in\Omega^{*}(B)$. 

The following projection formula can be easily derived from the definition of the push-forward (cf. \cite[eq.(1.16)]{berline_getzler_vergne}): for every $\eta\in\Omega^*(M)$ and $\beta\in\Omega^*(B)$,
\begin{equation}\label{eq_projection_formula}
\Pi_*(\Pi^*\beta\wedge\eta)=\beta\wedge\Pi_* \eta.
\end{equation}

Let us return to our setting of an $n$-dimensional euclidean vector space $V$. Recall our notation for the sphere bundle over $V$, $SV=V \times S^{n-1}$. And let $\pi:SV \to V$ be the projection map. Then 
\begin{displaymath}
\alpha|_{(x,v)}(w):=\langle v,d\pi_{(x,v)}w\rangle, \quad (x,v)\in SV,\, w\in T_{(x,v)}SV,
\end{displaymath} 
defines a \emph{contact form} on $SV$. 

To every compact convex subset $K \subset V$ one can associate its normal cycle $\nc(K)$, which is an $(n-1)$-dimensional Legendrian cycle in $SV$. Consequently we can integrate any form $\omega \in \Omega^{n-1}(SV)$ over $\nc(K)$. We obtain
\begin{equation}\label{nc_vanishes}\int_{\nc(K)} \alpha \wedge \tau=0 \quad \text{ and } \quad  \int_{\nc(K)} d\tau=0 \quad \text{ for every } \tau \in \Omega^{n-2}(SV).
\end{equation}
These two equations further imply that $\int_{\nc(K)} d\alpha \wedge \eta=0$ for $\eta \in \Omega^{n-3}(SV)$.  We refer to \cite{zaehle86} for the construction of the normal cycle of a convex body and its main properties. 

Let $0 \leq p \leq n-1$. Consider the partial flag manifold $F(n,p+1)=\{(v,E): v \in S^{n-1}, E \in \Gr_{p+1}(V), v \in E\}$. Denote by $\Pi$ the projection map  
\begin{displaymath}
 \Pi:V \times F(n,p+1) \to V \times S^{n-1}, (x,v,E) \mapsto (x,v).
\end{displaymath}
Set $r:=p(n-p-1)$ to be the dimension of the fiber of $\Pi$ and set $s:=n-1+r$. 

\begin{definition}\label{def:smoothFlagMeas}
A flag area measure $\Phi \in \FlagArea^{(p)}$ is called \emph{smooth} if there exists a translation-invariant form $\tau \in \Omega^s(\V\times F(n,p+1))^{tr}$ such that for any $f\in C^{\infty}(F(n,p+1))$ and for any $K \in \mathcal{K}(\V)$,
\begin{displaymath}
	\int_{F(n,p+1)} f \ d\Phi(K,\cdot)=\int_{\nc(K)} \Pi_*(f\tau). 
\end{displaymath}
The space of smooth flag area measures is denoted by $\FlagArea^{(p),sm}$.
\end{definition}
Note that since smooth functions are dense in the space of continuous functions, a flag area measure $\Phi$ is determined by the integrals $\int_{F(n,p+1)}f \ d\Phi(K,\cdot)$ with smooth functions $f$.

To each translation-invariant $s$-form $\tau$ on $V \times F(n,p+1)$ we associate a smooth flag area measure $\Phi$ as described in Definition \ref{def:smoothFlagMeas}. We now determine the kernel of the map $\Omega^s(\V \times F(n,p+1))^{tr} \to \FlagArea^{(p),sm}: \tau \mapsto \Phi$.  

The space $\Omega^s(V \times F(n,p+1))^{tr}$ admits a filtration as follows. For $(x,v,E) \in V \times F(n,p+1), 0\leq j\leq s$, we define
\begin{align*}
 \mathfrak{F}^{s,j}_{x,v,E}:=\{
 &\phi \in \largewedge^s T^*_{(x,v,E)}(V \times F(n,p+1))\,:\, \\
 & \forall f_1,\ldots,f_j \in T_{(x,v,E)} \Pi^{-1}(\Pi(x,v,E)), \ \phi(f_1,\ldots,f_j,-)=0 \},
\end{align*}
and 
\begin{align*}
 \mathfrak{F}^{s,j}:=\{
 &\omega \in \Omega^s(V \times F(n,p+1))^{tr} :\, \\
 &\forall (x,v,E) \in V \times F(n,p+1), \ \omega|_{(x,v,E)} \in \mathfrak{F}^{s,j}_{x,v,E} \}.
\end{align*}
Then 
\begin{displaymath}
 \Omega^s(\V \times F(n,p+1))^{tr} = \mathfrak{F}^{s,r+1} \supset \mathfrak{F}^{s,r} \supset \ldots \supset \mathfrak{F}^{s,0}=\{0\}.
\end{displaymath}

In the proof of the following theorem we will work with smooth convex bodies. Let us recall some general and well-known facts concerning these bodies. For a smooth compact convex body $K$ with strictly positive curvature, the normal cycle of $K$ is the graph of $(\mathrm{id},\nu):\partial K \to SV$, where $\nu:\partial K \to S^{n-1}$ is the Gauss map. 

Fix $x_0 \in \partial K$. Set $U:=T_{x_0}\partial K$ and $W:=U \oplus U$. The shape operator, also called the Weingarten map, $d\nu_{x_0}:U=T_{x_0}\partial K \to T_{\nu(x_0)} S^{n-1} \cong T_{x_0}\partial K=U$ is self-adjoint (see e.g.~\cite{spivak3_99}). Hence there is an orthonormal basis $u_1,\ldots,u_{n-1}$ of $U$ (the principal curvature directions) and $\lambda_1,\ldots,\lambda_{n-1}>0$ (the principal curvatures) such that $T_{(x_0,\nu(x_0))} \nc(K)$ is spanned by the vectors $w_i:=(u_i,\lambda_i u_i) \in W$. 

Conversely, let $(x_0,v_0) \in SV$ be given. Set $U:=v_0^\perp$ and $W:=U \oplus U$ endowed with the symplectic form $\theta((u_1,u_2),(\tilde u_1,\tilde u_2))=\langle u_1,\tilde u_2\rangle-\langle \tilde u_1,u_2\rangle$, $u_i,\tilde u_i\in U$. Let $u_1,\ldots,u_{n-1}$ be an orthogonal basis of $U$ and $\lambda_1,\ldots,\lambda_{n-1}>0$. The subspace $L$ of $W$ spanned by $(u_i,\lambda_i u_i), i=1,\ldots,n-1$, is a Lagrangian space. Lagrangian spaces of this type will be called {\it strictly positive}. Under these conditions, there exists a smooth compact convex body $K$ with strictly positive curvature such that $x_0 \in \partial K, \nu(x_0)=v_0$ and $T_{(x_0,\nu(x_0))} \nc (K)=L$. For instance,
\begin{displaymath}
K:=\left\{(x_1,\ldots,x_n): -1 \leq x_n \leq -\frac12 \sum_{i=1}^{n-1} \lambda_i x_i^2 \right\}
\end{displaymath} 
with respect to the coordinate system given by $(x_0;u_1,\ldots,u_{n-1},v_0)$ has this property.

\begin{theorem}\label{thm_kernel}
 A form $\tau \in \Omega^s(V \times F(n,p+1))^{tr}$ induces the trivial flag area measure if and only if 
 \begin{displaymath}
  \tau \in \langle \Pi^*\alpha, \Pi^* d\alpha,\mathfrak{F}^{s,r}\rangle.
 \end{displaymath}
\end{theorem}

\proof
Let us first show that the displayed space is a subset of the kernel. Let $\tau=\Pi^*\alpha \wedge \rho$ for some form $\rho$. Using \eqref{eq_projection_formula}, we obtain
\begin{displaymath}
\Pi_*(f \tau)=\Pi_*(\Pi^*\alpha \wedge f\rho)=\alpha \wedge \Pi_*(f \rho). 
\end{displaymath}
Hence the integral of this form over the normal cycle vanishes by \eqref{nc_vanishes}. The same argument shows that the forms $\Pi^*d\alpha \wedge \rho$ belong to the kernel. If $\tau \in \mathfrak{F}^{s,r}$, then $\Pi_*(f \tau)=0$ by the definition of the push-forward.  

Let us prove the other inclusion. Suppose that $\tau$ induces the trivial flag area measure. 
Fix $(x_0,v_0,E_0) \in V \times F(n,p+1)$ and set $U:=v_0^\perp$, which is an $(n-1)$-dimensional euclidean subspace of $V$.
Take a sequence of functions $h_j \in C^\infty(F(n,p+1))$ with $h_j(v_0,E_0) =1, h_j \geq 0$, whose supports shrink to $\{(v_0,E_0)\}$.
Then for all $K$, we have 
\begin{equation} \label{eq_kernel_condition}
 \int_{\nc(K)} \Pi_*(h_j \tau)=0. 
\end{equation}
Let us fix a strictly positive Lagrangian subspace $L$ of $W:=U \oplus U$ spanned by $w_i:=(u_i,\lambda_i u_i)$ as above. As noted above, there exists a smooth convex body $K$ with strictly positive curvature having $x_0 \in \partial K, \nu(x_0)=v_0$ and $T_{(x_0,v_0)} \nc(K)=L$. Let 
\begin{displaymath}
 A:=d\Pi|_{(x_0,v_0,E_0)}:T_{(x_0,v_0,E_0)} (\V \times F(n,p+1)) \to T_{(x_0,v_0)} (V \times S^{n-1}).
\end{displaymath} 

\noindent
Let $W_1,\ldots,W_{n-1}$ be lifts of $w_1,\ldots,w_{n-1}$, i.e. $W_i \in T_{(x_0,v_0,E_0)}(V\times F(n,p+1))$ and $A(W_i)=w_i, i=1,\ldots,n-1$. 
 Then \eqref{eq_kernel_condition}, the definition of the push-forward and a continuity argument imply that 
\begin{equation} \label{eq_kernel_condition2} 
\tau|_{(x_0,v_0,E_0)}(W_1,\ldots,W_{n-1},f_1,\ldots,f_r)=0,  \quad \forall f_1,\ldots,f_r \in \ker A.
\end{equation}

We may decompose  
\begin{displaymath}
 \tau|_{(x_0,v_0,E_0)} \equiv A^* \phi \wedge \kappa \mod \mathfrak{F}^{s,r}_{x_0,v_0,E_0},
\end{displaymath}
where $0 \neq \kappa \in \largewedge^r (\ker A)^*$ and $\phi \in \largewedge^{n-1} T^*_{(x_0,v_0)} (\V \times S^{n-1})$. For instance, we may use a local trivialization of $\Pi$, decompose $\tau|_{(x_0,v_0,E_0)}$ according to bidegrees and pick the highest degree part. Then \eqref{eq_kernel_condition2} implies that
\begin{displaymath}
 \phi(w_1,\ldots,w_{n-1})=0,
\end{displaymath}
i.e. $\phi$ vanishes on all strictly positive Lagrangian subsets of $W$.

Now we apply a technical Lemma \ref{lemma_bernig_broecker} below to the symplectic vector space $(W=U \oplus U,\theta:=-d\alpha_{(x_0,v_0)}|_W)$. Note that $T_{(x_0,v_0)}(V \times S^{n-1})=\R v_0 \oplus v_0^\perp \oplus v_0^\perp=\R v_0 \oplus W$ and that the projection onto $\R v_0$ equals $\alpha|_{(x_0,v_0)}$. Hence $\phi \in \langle \alpha|_{(x_0,v_0)},d\alpha|_{(x_0,v_0)}\rangle$. It follows that 
\begin{displaymath}
\tau|_{(x_0,v_0,E_0)} \in \langle A^*\alpha|_{(x_0,v_0,E_0)},A^*d\alpha|_{(x_0,v_0,E_0)},\mathfrak{F}^{s,r}_{x_0,v_0,E_0}\rangle.
\end{displaymath}
Since $x_0,v_0,E_0$ were arbitrary, this implies that $\tau \in \langle \Pi^*\alpha,\Pi^* d\alpha, \mathfrak F^{s,r}\rangle$, finishing the proof. 
\endproof

We complete the proof by an algebraic lemma which is adapted from \cite[Lemma 1.4]{bernig_broecker07}. Let $U$ be a euclidean vector space of dimension $m$ and let $W:=U \oplus U$ be endowed with its natural symplectic form $\theta$, $\theta((u_1,u_2),(\tilde u_1,\tilde u_2))=\langle u_1,\tilde u_2\rangle-\langle \tilde u_1,u_2\rangle$, $u_i,\tilde u_i\in U$. 

A subspace of the form $L:=\mathrm{span}\{w_1,\ldots,w_k\} \subset W$ with $w_i:=(u_i,\lambda_i u_i)$, $i=1,\dots,k$, where $u_1,\ldots,u_k$ are non-zero and pairwise orthogonal in $U$ and $\lambda_i>0$  is called a \emph{strictly positive isotropic subspace}.

\begin{lemma} \label{lemma_bernig_broecker}
Let $0 \leq k \leq 2m$ and suppose that $\omega \in \Lambda^k W^*$ vanishes on all strictly positive isotropic subspaces of dimension $k$. Then $\omega$ is a multiple of the symplectic form $\theta$.
\end{lemma}

\proof
We use induction on the dimension $m$ of $U$. The base case $m=1$ is trivial. 

Suppose $k>m$. Then basic symplectic algebra implies that every form $\omega \in \Lambda^k W^*$ is a multiple of $\theta$ and there is nothing to prove. 

Assume $k \leq m$. Let $u_1,\ldots,u_m$ be an orthonormal basis of $U$. Define a basis of $W$ by $e_i:=(u_i,0),f_i:=(0,u_i)$. Let $e_1^*,\ldots,e_m^*,f_1^*,\ldots,f_m^* \in W^*$ be the dual basis. Then the natural symplectic form on $W$ is $\theta=\sum_{i=1}^m e_i^* \wedge f_i^*$. Let $U_0$ be the euclidean subspace spanned by $u_1,\ldots,u_{m-1}$ and $W_0:=U_0 \oplus U_0$ with its natural symplectic form $\theta_0$. Then $\theta=\theta_0 + e_m^* \wedge f_m^*$. We decompose 
\begin{displaymath}
\omega=\omega_0 + \omega_1 \wedge e_m^* + \omega_2 \wedge f_m^* + \omega_3 \wedge e_m^* \wedge f_m^*,
\end{displaymath}
with $\omega_0 \in \Lambda^k W_0^*, \omega_1,\omega_2 \in \Lambda^{k-1} W_0^*, \omega_3 \in \Lambda^{k-2}W_0^*$.
Let $\lambda_1,\ldots,\lambda_{k-1},\lambda_m>0$ and let $v_1,\ldots,v_{k-1} \in U_0 \setminus \{0\}$ be pairwise orthogonal. Then $v_1,\ldots,v_{k-1},u_m$ are pairwise orthogonal in $U$ and the vectors $w_i:=(v_i,\lambda_i v_i), i=1,\ldots,k-1$ and $w_m:=(u_m,\lambda_m u_m)$ span a strictly positive isotropic space. By assumption, 
\begin{displaymath}
0=\omega(w_1,\ldots,w_{k-1},w_m)=\omega_1(w_1,\ldots,w_{k-1})+\lambda_m \omega_2(w_1,\ldots,w_{k-1}).
\end{displaymath} 
With $\lambda_m>0$ being arbitrary, it follows that 
\begin{displaymath}
\omega_1(w_1,\ldots,w_{k-1})=0,\quad \omega_2(w_1,\ldots,w_{k-1})=0.
\end{displaymath}
By the induction hypothesis, $\omega_1,\omega_2$ are multiples of $\theta_0$. Since $\theta_0 \wedge e_m^{*}=\theta \wedge e_m^*$, $\theta_0 \wedge f_m^* = \theta \wedge f_m^*$, we obtain that $\omega_1 \wedge e_m^*+\omega_2 \wedge f_m^*$ is a multiple of $\theta$.

Let $v_1,\ldots,v_k \in U_0 \setminus \{0\}$ be pairwise orthogonal unit vectors and $w_i:=(v_i,\lambda_i v_i)$ with $\lambda_i>0$. Then $0=\omega_0(w_1,\ldots,w_k)$ and the induction hypothesis implies that $\omega_0=\omega_{00} \wedge \theta_0$ with $\omega_{00} \in \Lambda^{k-2}W_0^*$. Let $\lambda>0$ and $\tilde w_{k-1}:=(v_{k-1}+u_m,\lambda(v_{k-1}+u_m)) \in W, \tilde w_k:=(v_{k-1}-u_m, v_{k-1}-u_m) \in W$. Then, since $\omega_1\wedge e_m^*+\omega_2\wedge f_m^*$ is a multiple of $\theta$, 
\begin{align*}
0 & = (\omega_0+\omega_3 \wedge e_m^* \wedge f_m^*)(w_1,\ldots,w_{k-2},\tilde w_{k-1},\tilde w_k) \\
& =  (\lambda-1) (\omega_3-\omega_{00})(w_1,\ldots,w_{k-2}).
\end{align*}
By the induction hypothesis $\omega_3-\omega_{00}=\theta_0 \wedge \phi_0$ with $\phi_0 \in \Lambda^{k-4} W_0^*$. Hence  
\begin{align*}
\omega_0+\omega_3 \wedge e_m^* \wedge f_m^* & =\omega_{00} \wedge \theta_0+(\omega_{00}+\theta_0 \wedge \phi_0) \wedge e_m^* \wedge f_m^* \\
& = (\omega_{00} +\phi_0 \wedge e_m^* \wedge f_m^*) \wedge \theta.
\end{align*} 
This shows that $\omega_0+\omega_3 \wedge e_m^* \wedge f_m^*$ is divisible by $\theta$.
\endproof

\section{Dimension computation}\label{sec_dimension}

In this section we determine the dimension of the space of smooth, $k$-homogeneous, translation- and $\SO(n)$-covariant valuations with values in the space of signed measures on the partial flag manifold $F(n,p+1)$, denoted by $\FlagArea^{(p),sm,\mathrm{SO}(n)}_k$.

\begin{theorem} Let $0 \leq p, k \leq n-1$. Then 
\begin{align*}
 \dim \FlagArea^{(p),sm,\mathrm{SO}(n)}_k & = \min\{k,n-k-1,p,n-p-1\}+1 +\\
 & \quad + \begin{cases} 1 & \text{ if }p=k=\frac{n-1}{2},\\ 0 & \text{otherwise}. \end{cases}
\end{align*}
\end{theorem}

\proof
If $p=0$ or $p=n-1$, then $F(n,p+1) \cong S^{n-1}$ and the result follows from Schneider's characterization of area measures \cite{schneider78}. For the rest of the proof, we assume that $1 \leq p \leq n-2$. 

Since $F(n,p+1) \cong F(n,q+1)$ with $q:=n-p-1$, there is an isomorphism 
\begin{displaymath} 
\FlagArea_k^{(p),sm,\overline{\SO(n)}} \cong \FlagArea_k^{(q),sm,\overline{\SO(n)}}.
\end{displaymath}
We may thus assume that $1 \leq p \leq q$.

According to Definition \ref{def:smoothFlagMeas}, a smooth flag area measure is obtained by integration over the normal cycle of the push-forward of a form $\tau$ on $V \times F(n,p+1)$. Since all our forms will be obtained as a push-forward of appropriate forms on $\overline{\SO(n)} = V \rtimes \SO(n)$, let us describe the different spaces on which we construct the forms and the maps between them.

Let $G:=\SO(n), \overline G:=\overline{\SO(n)}$, and $H:=S(\O(p) \times \O(q))$. Then $F(n,p+1) \cong G/H$ and $V \times F(n,p+1) \cong \overline{G}/H$. Fix an orthonormal basis $e_1,\dots,e_n$ of $V$. For an element $g \in G$, we let $g_i=ge_i$, $1\leq i\leq n$, denote the column vectors of $g$. Let 
\begin{displaymath}
 \tilde \Pi: V \rtimes \SO(n) \to V \times F(n,p+1)
\end{displaymath}
be the projection map, taking $(x,g)\mapsto (x,g_1, \spann\{g_1, \ldots, g_{p+1} \})$. The fiber of this map is diffeomorph to $H$ and has dimension $\tilde{r}:=\frac{p(p-1)+q(q-1)}{2}$. 

Let 
\begin{displaymath}
 \Pi: V \times F(n,p+1) \to V \times S^{n-1}
\end{displaymath}
be the projection map, taking $(x, v,E) \mapsto (x,v)$. The fiber of this map is diffeomorph to $\Gr_p(\R^{n-1})$ and has dimension $r:=p(n-1-p)=pq$. 

Finally, let 
\begin{displaymath}
 \hat{\Pi}: V \rtimes \SO(n) \to V \times S^{n-1}
\end{displaymath}
be the projection map, taking $(x, g) \mapsto (x,g_1)$. The fiber of this map is diffeomorph to $\SO(n-1)$ and has dimension $\hat{r}:=\frac{(n-1)(n-2)}{2}$. Thus we obtain the following commutative diagram:
\begin{displaymath}
 \xymatrix{\overline G = V \rtimes \SO(n) \ar[r]^-{\tilde \Pi} \ar[dr]_-{\hat \Pi} & V \times F(n,p+1) \ar[d]^{\Pi} \\ & V \times S^{n-1}}
\end{displaymath}

A basis of the algebra of left-invariant forms on $\overline G$ is given by the $1$-forms $\sigma_i, 1\leq i \leq n$ and $\omega_{ij}=-\omega_{ji}, 1 \leq i < j \leq n$, see Subsection \ref{subsec_construction_forms} for their definition. These $1$-forms span the dual space $\overline{\mathfrak{g}}^*$. Denote by $X_i, 1\leq i \leq n, X_{ij}, 1 \leq i<j \leq n$ the corresponding basis of the Lie algebra $\overline{\mathfrak{g}}$. 

Denote the linear span of $X_1$ by $V_0$; the linear span of $X_i$ with $2 \leq i \leq p+1$ by $V_\sigma$ and the linear span of $X_i$ with $p+2 \leq i \leq n$ by $W_\sigma$. The linear span of $X_{1j}$ with $2 \leq j \leq p+1$ will be denoted by $V_\omega$; the linear span of $X_{1j}$ with $p+2 \leq j \leq n$ by $W_\omega$; the linear span of $X_{ij}$ with $2 \leq i<j\leq p+1$ by $U_p$; the linear span of $X_{ij}$ with $2 \leq i < p+2 \leq j\leq n$ by $U_{pq}$ and the linear span of $X_{ij}$ with $p+2 \leq i<j\leq n$ by $U_q$. Note that $\ker d\tilde \Pi=U_p \oplus U_q$, $\ker d \Pi=U_{pq}$ and $\ker d\hat \Pi=U_p \oplus U_q \oplus U_{pq}$. Schematically, the Lie algebra $\overline{\mathfrak{g}}$ looks as follows:
\begin{displaymath}
 \overline{\mathfrak{g}}=\left(\begin{matrix} 0 & & & \\ V_0 & 0 & & \\V_\sigma & V_\omega & U_p & \\ W_\sigma & W_\omega & U_{pq} & U_q\end{matrix}\right).
\end{displaymath}

A smooth flag area measure is generated from a translation-invariant form $\omega \in \Omega^{n-1+r}(\overline{G}/H)^{tr}$. A $k$-homogeneous smooth flag area measure corresponds to a form of bidegree $(k,n-1+r-k)$, $\omega \in \Omega^{k,n-1+r-k}(\overline{G}/H)^{tr}$, where the bidegree is taken with respect to the product structure of $\overline{G}/H=V \times G/H$. For a $G$-covariant flag area measure we may assume by averaging with respect to the Haar measure on the compact group $G$ that the form $\omega$ is $G$-invariant. 

Thus let $\omega \in \Omega^{k,n-1+r-k}(\overline{G}/H)^{\overline G}$ and set $\tilde \omega:=\tilde \Pi^*\omega \in \Omega^{k,n-1+r-k}(\overline G)$. The form $\tilde\omega$ is left $\overline G$-invariant and right $H$-invariant. It vanishes on each tangent vector to the fiber of $\tilde \Pi$, i.e. on $U_p \oplus U_q$. Conversely, every left $\overline G$-invariant and right $H$-invariant form on $\overline G$, which vanishes on $U_p \oplus U_q$, induces a $\overline G$-invariant form on $\overline{G}/H$. 

Since $(\tilde \Pi^* \circ \Pi^*)\alpha=\hat \Pi^*\alpha=\sigma_1$, the quotient of the space of $\overline G$-invariant forms on $\overline{G}/H$ by multiples of $\Pi^* \alpha$ is identified with the space 
\begin{multline*}
 (\largewedge^k (V^*_\sigma \oplus W^*_\sigma) \otimes \largewedge^{n-1+r-k} (V^*_\omega \oplus W^*_\omega \oplus U^*_{pq}))^H \\
 = \bigoplus_{i=r-k}^{r} \left[\largewedge^k 
  (V^*_\sigma \oplus W^*_\sigma) \otimes \largewedge^{n-1+r-k-i}(V^*_\omega \oplus W^*_\omega) \otimes \largewedge^i U^*_{pq}\right]^H.
\end{multline*}

If $\tilde \omega$ belongs to the sum of terms with $i<pq=r$, then $\omega \in \mathfrak F^{n-1+r,r}$ and hence $\tilde \omega$ induces the trivial flag area measure by Theorem \ref{thm_kernel}. Moreover, $\largewedge^{pq}U_{pq}^*$ is the trivial one-dimensional representation of $H$. We obtain 
\begin{displaymath}
	\Omega^{k,n-1+r-k}(\overline{G}/H)^{\overline G}/\langle \Pi^*\alpha,\mathfrak F^{n-1+r,r}\rangle \cong \left[\largewedge^k 
  (V^*_\sigma \oplus W^*_\sigma) \otimes \largewedge^{n-1-k}(V^*_\omega \oplus W^*_\omega)\right]^H.
\end{displaymath}

We also have to quotient out the form $\Pi^*d\alpha$. Note that $(\tilde \Pi^* \circ \Pi^*) d\alpha=\hat \Pi^* d\alpha=-\sum_{i=2}^n \sigma_i \wedge \omega_{1i}$. This is a symplectic form on the $2(n-1)$-dimensional vector space $(V_\sigma \oplus W_\sigma) \oplus (V_\omega \oplus W_\omega)$. 
By basic symplectic linear algebra (cf. \cite[Prop. 1.2.30]{huybrechts05}), multiplication by this form induces an injective map 
\begin{multline*}
	L:\underbrace{\left[\largewedge^{k-1} 
  (V^*_\sigma \oplus W^*_\sigma) \otimes \largewedge^{n-2-k}(V^*_\omega \oplus W^*_\omega)\right]^H}_{=:B_{k,p,q}} \\
  \to \underbrace{\left[\largewedge^k 
  (V^*_\sigma \oplus W^*_\sigma) \otimes \largewedge^{n-1-k}(V^*_\omega \oplus W^*_\omega)\right]^H}_{=:A_{k,p,q}}.
\end{multline*}
Therefore we may view $B_{k,p,q}$ as a subspace of $A_{k,p,q}$ and obtain 
\begin{displaymath}
 \Omega^{k,n-1+r-k}(\overline{G}/H)^{\overline G}/\langle \Pi^*\alpha,\Pi^*d\alpha,\mathfrak F^{n-1+r,r}\rangle \cong A_{k,p,q}/B_{k,p,q}.
\end{displaymath}
By Theorem \ref{thm_kernel}, it follows that 
\begin{align*}
 \dim \FlagArea^{(p),sm,\mathrm{SO}(n)}_k & = \dim \Omega^{k,n-1+r-k}(\overline{G}/H)^{\overline G}/\langle \Pi^*\alpha,\Pi^*d\alpha,\mathfrak F^{n-1+r,r}\rangle\\
 & =a_{k,p,q}-b_{k,p,q},
\end{align*}
where  
\begin{align*}
a_{k,p,q} & := \dim A_{k,p,q} \\
& = \dim \left[ \largewedge^k(V^*_{\sigma} \oplus W^*_{\sigma}) \otimes  \largewedge^{n-k-1}(V^*_{\omega} \oplus W^*_{\omega}) \right]^H \\
& = \sum_{i,j} \dim \left[\largewedge^i V^*_{\sigma} \otimes \largewedge^{k-i} W^*_{\sigma} \otimes \largewedge^j V^*_{\omega} \otimes \largewedge^{n-k-1-j} W^*_{\omega} \right]^H
\end{align*}
with $\max\{0, k-q\} \leq i \leq \min\{k,p\}, \max\{0, p-k\} \leq j \leq \min\{n-k-1,p\}$ and 
\begin{align*}
b_{k,p,q} & := \dim B_{k,p,q} \\
& = \dim \left[ \largewedge^{k-1}(V^*_{\sigma} \oplus W^*_{\sigma}) \otimes  \largewedge^{n-k-2}(V^*_{\omega} \oplus W^*_{\omega}) \right]^H  \\
& =\sum_{i,j} \dim \left[\largewedge^i V^*_{\sigma} \otimes \largewedge^{k-1-i} W^*_{\sigma} \otimes \largewedge^j V^*_{\omega} \otimes \largewedge^{n-k-2-j} W^*_{\omega} \right]^H
\end{align*}
with $\max\{0, k-1-q\} \leq i \leq \min\{k-1,p\}, \max\{0, p-k-1\} \leq j \leq \min\{n-k-2,p\}$. 

Let us first consider the action of the subgroup $H':=\mathrm{SO}(p) \times \mathrm{SO}(q) \subset H$. We have 
\begin{align} 
&\left[\largewedge^i V^*_{\sigma} \otimes \largewedge^{k-i} W^*_{\sigma} \otimes \largewedge^j V^*_{\omega} \otimes \largewedge^{n-k-1-j} W^*_{\omega} \right]^{H'} \label{eq_sop_soq} \\
&=\left[\largewedge^i V^*_{\sigma} \otimes \largewedge^j V^*_{\omega} \right]^{\SO(p)} \otimes \left[ \largewedge^{k-i} W^*_{\sigma} \otimes \largewedge^{n-k-1-j} W^*_{\omega} \right]^{\SO(q)}. \nonumber
\end{align}
We decompose 
\begin{displaymath}
(\largewedge^i V^*_\sigma \otimes \largewedge^j V^*_\omega)^{\SO(p)} = \bigoplus_{\epsilon=0,1} (\largewedge^i V^*_\sigma \otimes \largewedge^j V^*_\omega)^{\SO(p),\epsilon},
\end{displaymath}
where $(\largewedge^i V^*_\sigma \otimes \largewedge^j V^*_\omega)^{\SO(p),\epsilon}$ is the space of $\SO(p)$-invariant elements such that $g \in \O(p)$ acts by $\det(g)^\epsilon$. An easy exercise in representation theory (see \cite[Lemma 5.1]{bernig_faifman_opq} and \cite[Lemma 0.4.3]{fu90}) is to show that  
\begin{align*}
(\largewedge^i V^*_\sigma \otimes \largewedge^j V^*_\omega)^{\SO(p),+} & = \delta_{i}^j,\\
(\largewedge^i V^*_\sigma \otimes \largewedge^j V^*_\omega)^{\SO(p),-} & = \delta_{i+j}^p.
\end{align*}
More precisely, for $i=j$, the space $(\largewedge^i V^*\otimes \largewedge^i V^*)^{\SO(p),+}$ is spanned by the $i$-th power of the symplectic form on $V$ and for $i+j=p$, $(\largewedge^i V^* \otimes \largewedge^{p-i} V^*)^{\SO(p),-}$ is spanned by the determinant.

The space $\left[ \largewedge^{k-i} W^*_{\sigma} \otimes \largewedge^{n-k-1-j} W^*_{\omega} \right]^{\SO(q)}$ can be decomposed in an analogous way. An $H'$-invariant element is invariant under the larger group $H$ if and only if the factors in the decomposition \eqref{eq_sop_soq} have the same $\epsilon$. Therefore  
\begin{displaymath}
 a_{k,p,q} = \sum_{i=\max\{0, k-q\}}^{\min\{k,p\}} \sum_{j=\max\{0, p-k\}}^{\min\{n-k-1,p\}} (\delta_i^j \delta_{k-i}^{n-k-1-j}+\delta_{i+j}^p \underbrace{\delta_{k-i+n-k-1-j}^q}_{=\delta_{i+j}^p}).
\end{displaymath}

The first summand is non-zero only if $i=j$ and $k-i=n-k-1-j$, i.e. if $i=j$ and $k=n-k-1$. In this case, both sums range from $0$ to $p$. By treating similarly the second summand, we have
\begin{displaymath}
a_{k,p,q}  = \delta_k^{n-k-1} (p+1) + \min\{k,n-1-k,p,q\}+1.
 \end{displaymath}

Similar arguments yield 
\begin{align*}
 b_{k,p,q}  & = \sum_{i=\max\{0, k-1-q\}}^{\min\{k-1,p\}} \sum_{j=\max\{0, p-k-1\}}^{\min\{n-k-2,p\}}  \delta_i^j \delta_{k-1-i}^{n-k-2-j} \\
 & = \delta_{k}^{n-k-1} \min\left\{\frac{n-1}{2},p+1\right\}.
\end{align*}

Overall we obtain that 
\begin{align*}
 \dim &\FlagArea^{(p),sm,\mathrm{SO}(n)}_k  =  a_{k,p,q}-b_{k,p,q}\\
 & =\min\{k,n-k-1,p,q\}+1+\delta_k^{n-k-1} \left[(p+1)-\min\left\{\frac{n-1}{2},p+1\right\}\right]\\
 & = \min\{k,n-k-1,p,q\}+1+\begin{cases} 1 & p=k=\frac{n-1}{2},\\ 0 & \text{otherwise}, \end{cases}
\end{align*}
as claimed.
\endproof

\section{Jordan angles}\label{sec_jordan}

In this section we introduce the notion of angles between subspaces and discuss their properties. Let $\W$ be a euclidean vector space of dimension $n'$ and let $E$ and $F$ be subspaces of $\W$ of dimensions $p$ and $k$, respectively. We start by recalling the notion of the cosine between subspaces. Let us denote by $B_{E,F}$ the orthogonal projection from $E$ onto $F$. In case $p \leq k$, let $S$ be a domain in $E$ of volume $1$. The absolute value of the cosine between the subspaces $E,F$ is defined as the $p$-dimensional volume of the orthogonal projection of $S$ onto $F$:
$$ |\cos(E,F)|:= \vol_p(B_{E,F} S) .$$
In case $p>k$, we define the cosine between $E$ and $F$ analogously, by considering the orthogonal projection from $F$ onto $E$. 

The cosine is invariant under the diagonal action of $\O(n')$ on the product of two Grassmannians, but it is not enough to separate the orbits. For that we use a more general notion of {\it Jordan angles} between subspaces, sometimes also called {\it critical angles} or {\it principal angles}, see for e.g. \cite{GallagherProulx_1977,jordan1875}. To define Jordan angles associated with a pair of subspaces $(E,F)$, we need to start with appropriate bases in our subspaces. Those bases are given by the following simple lemma.

\begin{lemma}[{\cite[Lemma 1]{GallagherProulx_1977}}] \label{lem_orbits_grassmann}
Let $E$ and $F$ be subspaces of an $n'$-dimensional euclidean space $\W$ of dimensions $p$ and $k$, respectively. Set $m:=\min\{k,n'-k,p,n'-p\}$. Then there is an orthonormal basis $e_1, \ldots, e_p$ of $E$, an orthonormal basis $f_1, \ldots, f_k$ of $F$ and $\frac{\pi}{2} \geq \theta_1 \geq \ldots \geq \theta_m \geq 0$ such that 
\begin{enumerate}
  \item $\left\langle  e_i, f_j \right\rangle=0$ for $i \neq j$.
   \item $\left \langle e_i,f_i\right\rangle=\cos \theta_i$ for $1 \leq i \leq m$.
 \item $\left\langle  e_i, f_i \right\rangle=1$ for $m+1 \leq i \leq \min\{k,p\}$.
\end{enumerate}
The unique angles $\theta_1,\ldots,\theta_m$ are called Jordan angles between the subspaces $E$ and $F$.
 
Two pairs $(E,F),(E',F') \in \Gr_p(W) \times \Gr_k(W)$ belong to the same $\mathrm{O}(n')$-orbit if and only if the Jordan angles between $E$ and $F$ and between $E'$ and $F'$ are the same. 
\end{lemma}

In the following lemma we construct convenient bases for the subspaces $E, F, E^{\perp}$ and $F^{\perp}$ simultaneously.

\begin{lemma} \label{lem_orbits_grassmann_orthogonal}
With $E, F, \{e_i\}, \{f_j\}$ as in the previous lemma, there are orthonormal bases $h_1,\ldots,h_{n'-k}$ of $F^\perp$ and $g_1,\ldots,g_{n'-p}$ of $E^\perp$ such that 
\begin{enumerate}
 \item $\left \langle e_i,h_i\right\rangle= \left\langle g_i,f_i\right\rangle =\sin\theta_i , \left \langle g_i,h_i\right\rangle=-\cos\theta_i$ for $1 \leq i \leq m$.
 \item $g_i=h_i, m+1 \leq i \leq \min\{n'-k,n'-p\}$.
 \item If $p \geq k$, then $h_{n'-p+i}=e_{k+i}, 1 \leq i \leq p-k$.
 \item If $k \geq p$, then $g_{n'-k+i}=f_{p+i}, 1 \leq i \leq k-p$.
 \item All other scalar products are zero.  
\end{enumerate}
In particular, the Jordan angles between $E$ and $F^\perp$ are given by $\theta_j'=\frac{\pi}{2}-\theta_{m-j+1}, j=1,\ldots,m$.
\end{lemma}

\proof
Let $0 \leq l \leq m$ be the largest index with $\theta_l \neq 0$. Then a basis of $E \cup F$ is given by 
\begin{displaymath}
\{e_i:i=1,\ldots,l\} \cup \{f_i:i=1,\ldots,l\} \cup \{e_i= f_i: i=l+1,\ldots,\min\{k,p\}\}
\end{displaymath}
together with 
\begin{displaymath}
\{e_i:k+1 \leq i \leq p\} \text{ if } p > k \quad \quad \text{ or } \quad \quad \{f_i:p+1 \leq i \leq k\} \text{ if } p < k.
\end{displaymath}
Hence $\dim(E\cup F)=\max\{k,p\}+l$. 

For $1 \leq i \leq l$, let $h_i$ be the unit vector in the oriented $2$-dimensional plane spanned by the vectors $e_i,f_i$ obtained from $f_i$ by a rotation by $-\frac{\pi}{2}$ and let $g_i$ be the unit vector obtained by rotating $e_i$ by $\frac{\pi}{2}$. Then $\langle h_i,e_i\rangle= \langle f_i,g_i\rangle=\sin \theta_i$ and $ \langle h_i,g_i\rangle=-\cos \theta_i$. 

For $p \leq k$, $\dim(E^\perp \cap F^\perp)=\dim(E\cup F)^\perp = n'- \dim(E\cup F) = n'-k-l$. Let $h_i, i=l+1,\ldots,n'-k$ be an orthonormal basis of $E^\perp \cap F^\perp$. Vectors $\{h_i\}$ form a basis of $ F^\perp$. To complete a basis of $E^\perp$, set $g_i:=-h_i, i=l+1,\ldots, m$, $g_i:=h_i, i=m+1,\ldots,n'-k$ and $g_{n'-k+i}:=f_{p+i}, i=1,\ldots,k-p$. 

Similarly, for $p \geq k$, $\dim(E^\perp \cap F^\perp) = n'-p-l$. Let $g_i, i=l+1,\ldots,n'-p$ be an orthonormal basis of $E^\perp \cap F^\perp$. Vectors $\{g_i\}$ form a basis of $ E^\perp$. To complete a basis of $F^\perp$, set $h_i:=-g_i, i=l+1,\ldots, m$, $h_i:=g_i, i=m+1,\ldots,n'-p$ and $h_{n'-p+i}:=e_{k+i}, i=1,\ldots,p-k$. 
\endproof

\begin{definition}\label{def_sigmai}
Let $E \in \Gr_p(W), F \in \Gr_k(W)$ have Jordan angles $\theta_1,\ldots,\theta_m$, where $m=\min\{k,n'-k,p,n'-p\}$. We define $\sigma_i(E,F)$ to be the $i$-th elementary symmetric function in $\cos^2 \theta_1,\ldots,\cos^2 \theta_m$, for $0 \leq i \leq m$.
\end{definition}

The following properties are obvious. 
\begin{enumerate}
\item $\sigma_i$ is $\O(n')$-invariant, i.e. $\sigma_i(g E, g F) = \sigma_i(E, F)$ for all $g \in \O(n')$.
\item $\sigma_i(E,F)=\sigma_i(F,E)$. 
\item $\sigma_m(E,F)$ equals the squared cosine between $E$ and $F$; and $\sigma_m(E^\perp,F)$ equals the squared sine between $E$ and $F$.
\end{enumerate}

The next result gives the probability distribution of the Jordan angles between a fixed plane and a random plane in $\W$.

\begin{theorem}[{\cite[Section 6]{james54}}] \label{thm_james}
Let $\W$ be an $n'$-dimensional euclidean vector space. Let $F\in\Gr_k(\W)$ be fixed and let $E \in \Gr_p(\W)$ be chosen randomly according to the $\SO(n')$-invariant probability measure. Denote the Jordan angles between $E$ and $F$ by $\theta_i$, $1\leq i \leq m$, and set $x_i:=\cos^{2}(\theta_i)$. Then the probability density of the $(x_1,\ldots,x_m)$ is proportional to  
\begin{displaymath}
\prod_{j=1}^m x_j^{\frac{|p-k|-1}{2}}(1-x_j)^{\frac{|n'-p-k|-1}{2}}\prod_{1\leq i<j\leq m}(x_j-x_i).
\end{displaymath}
\end{theorem}

In \cite[Section 6]{james54}, this formula is shown under the assumption $k \leq p \leq \frac{n'}{2}$, but using Lemma \ref{lem_orbits_grassmann_orthogonal}, the other cases can be checked as well.

We note that the integral of a symmetric function in $x_1,\ldots,x_m$ with respect to this density over $\{0\leq x_1\leq\dots\leq x_m\leq 1\}$ can be written as an integral over $[0,1]^m$ by dividing the integral by $m!$ and replacing the factor $\prod_{1\leq i<j\leq m}(x_j-x_i)$ by $\prod_{1\leq i<j\leq m}|x_j-x_i|$. Such integrals were studied by Selberg \cite{selberg44} and others, see the survey \cite{forrester_warnaar}. In this paper, we will need the following integral of Selberg type. 

\begin{theorem}[{\cite[Theorem 2]{aomoto87}}] \label{thm_aomoto}
Let $m \in \N$, $\lambda>0$, $\lambda',\lambda''>-1$ and let $f:\mathbb{R}^m \to \mathbb{R}$. Define 
\begin{displaymath}
J_f=\int_{[0,1]^m} f(x_1,\dots,x_m) \prod_{j=1}^m x_j^{\lambda'}(1-x_j)^{\lambda''} \prod_{1\leq i<j\leq m} |x_i-x_j|^{\lambda}dx_1\dots dx_m.
\end{displaymath}
If $f(x_1,\dots,x_m):=\prod_{j=1}^m(x_j-t)$ for $t \in \mathbb{R}$, then 
\begin{equation}\label{aomoto}
\frac{J_f}{J_1}=\sum_{r=0}^m(-t)^{m-r}\binom{m}{r} \prod_{j=1}^{r}\frac{\lambda'+1+\frac12(m-j)\lambda}{\lambda'+\lambda''+2+\lambda(m-\frac{j}{2}-\frac12)}.
\end{equation} 
\end{theorem}

\begin{corollary} \label{cor_expectation_sigma_i}
Let $\W$ be an $n'$-dimensional euclidean vector space. Let $F \in \Gr_k(\W)$ be fixed and let $E \in \Gr_p(\W)$ be chosen randomly according to the $\SO(n')$-invariant probability measure. Then the expectation of $\sigma_i(E,F), 0 \leq i \leq m$ is given by 
\begin{displaymath}
 \binom{m}{i} \binom{|p-k|+m}{i} \binom{n'}{i}^{-1}.
\end{displaymath}
\end{corollary}

\proof
By Theorems \ref{thm_james} and \ref{thm_aomoto}, we have to compute $\frac{J_{\sigma_i}}{J_1}$, where $\lambda=1$, $\lambda'=\frac{|p-k|-1}{2}$, $\lambda''=\frac{|n'-p-k|-1}{2}$. Since
\begin{equation}\label{pol2}
f(x_1,\dots,x_m)=\prod_{j=1}^m (x_j-t)=\sum_{r=0}^m(-t)^{m-r}\sigma_r(x_1,\dots,x_m),
\end{equation} 
Theorem \ref{thm_aomoto} implies that 
\begin{align*}
 \frac{J_{\sigma_i}}{J_1} & = \binom{m}{i} \prod_{j=1}^{i}\frac{\lambda'+1+\frac12(m-j)\lambda}{\lambda'+\lambda''+2+\lambda(m-\frac{j}{2}-\frac12)}\\
 & = \binom{m}{i} \binom{2\lambda'+m+1}{i} \binom{2\lambda'+2\lambda''+2m+2}{i}^{-1}\\
 & = \binom{m}{i} \binom{|p-k|+m}{i} \binom{|p-k|+|n'-p-k|+2m}{i}^{-1}\\
 & = \binom{m}{i} \binom{|p-k|+m}{i} \binom{n'}{i}^{-1},
\end{align*}
where the last line follows by 
\begin{align*}
 |p-k| & +|n'-p-k|  +2m \\
 & = \max\{p-k,k-p\}+\max\{n'-p-k,p+k-n'\}\\
 & \quad +2\min\{k,n'-k,p,n'-p\}\\
 & = \max\{n'-2k,2p-n',n'-2p,2k-n'\}+2\min\{k,n'-k,p,n'-p\}\\
 & = 2(\max\{k,n'-k,p,n'-p\}+\min\{k,n'-k,p,n'-p\})-n'\\
 & = n'.
\end{align*}
\endproof

We end this section with a definition of the angle between two subspaces both having dimension equal to the half of the dimension of the ambient euclidean space. 

\begin{definition}\label{def_cosinus}
Let $\W$ be an oriented euclidean vector space of dimension $n'=2a$ and let $E,F \in \Gr_a(\W)$.  Fix some orientations of $E$ and $F$, and endow $F^\perp$ with the orientation such that $W \cong F^\perp \oplus F$ is orientation preserving. 

Let $B_{E,F}$ (resp.~$B_{E,F^{\perp}}$) denote the orthogonal projection from $E$ onto $F$ (resp.~from $E$ onto $F^{\perp}$). We define 
\begin{displaymath}
\tilde\sigma_{a}(E,F)=\det(B_{E,F})\det(B_{E,F^{\perp}}), 
\end{displaymath}
where $\det(B_{E,F})$ denotes the determinant of the map $B_{E,F}$.
\end{definition}

It is easy to see that $\tilde \sigma_a(E,F)$ is independent of the choice of the orientations of $E$ and $F$, but changes its sign when we reverse the orientation of $W$.

Moreover, it is easy to check that
\begin{displaymath}
 |\tilde \sigma_a(E,F)|=\cos(E,F) \cos(E, F^\perp)
\end{displaymath}
and that $\tilde \sigma_a$ is invariant under the diagonal action of $\SO(2a)$. For $g\in\O(2a)$, we have
$$\tilde\sigma_a(gE,gF)=\det g\,\tilde\sigma_a(E,F).$$
Moreover, 
\begin{equation} \label{eq_tildesigma_perp}
\tilde\sigma_a(E,F)=(-1)^a \tilde\sigma_a(E^\perp,F),
\end{equation}
and
\begin{equation} \label{eq_tildesigma_comm}
\tilde\sigma_a(E,F)=(-1)^a \tilde\sigma_a(F,E).
\end{equation}
Both equations follow by using that the orientations of $E \oplus E^\perp$ and of $E^\perp \oplus E$ differ by a factor $(-1)^a$, and the fact that $|\det (B_{E,F})|=|\det(B_{E^{\perp},F^{\perp}})|$, since the block matrix with blocks $B_{E,F},B_{E^{\perp},F},B_{E,F^{\perp}}$, and $B_{E^{\perp},F^{\perp}}$ is a special orthogonal matrix.

\section{Construction of invariant flag area measures}\label{sec_construction}

The aim of this section is to prove the main results of this paper: Theorems~\ref{mainthm_construction}, \ref{mainthm_smooth}, \ref{mainthm_properties} and \ref{thm_basis}. The idea is to follow closely the proof of the dimension formula from Section \ref{sec_dimension} and to construct the invariant forms in an explicit way. 

\subsection{Construction of forms}\label{subsec_construction_forms}

Let $\overline G=\overline{\SO(n)}=\V \rtimes \SO(n)$ be the euclidean motion group and let $\pi_1:\overline G \to \V$, $\pi_2:\overline G \to \SO(n)$ be the projections onto the first and second factor. The Maurer-Cartan form on $\overline G$ takes values in the Lie algebra $\bar{\mathfrak g}= \V \rtimes \mathfrak{so}(n)$ (see e.g. \cite{spivak2_99}). 

Let us fix an orthonormal basis $e_1,\dots,e_n$ of $\V$. Let $\sigma_i$, $\omega_{ij}$ be the components of the Maurer-Cartan form with respect to this basis. Then $\omega_{ij}=-\omega_{ji}$. The forms $\omega_{ij}$, $1\leq i<j \leq n$, together with $\sigma_i$, $1\leq i\leq n$ form a basis of the space of left-invariant $1$-forms on $\overline G$.

For $\bar g \in \overline G$, let $\pi_2(\bar g)=:g=(g_1,\dots,g_n) \in G=\SO(n)$ be its rotational part, where $g_i=ge_i$, $1\leq i\leq n$, denotes the $i$-th column of $g$. We may consider $g_i$ as a vector valued function on $\overline G$. By definition (see \cite{santalo76}),
\begin{align*}
 \omega_{ij}|_{\bar g}(v) & =\langle g_i,dg_j(v)\rangle,\\
 \sigma_i|_{\bar g}(v) & = \langle g_i,d\pi_1(v)\rangle, \quad v \in T_{\bar g} \overline G.
\end{align*}

Let $1\leq k, p \leq n-1$. Set $q:=n-1-p$. The partial flag manifold $F(n,p+1)$ can be identified with the homogeneous space $G/H$, where $H=S(\O(p)\times\O(q))$, as described in Section~\ref{sec_dimension}. 

Let $\Pi,\tilde \Pi,\hat \Pi$ be the maps from Section \ref{sec_dimension}. We denote the  volume form on the corresponding fiber by $\rho,\tilde\rho$ and $\hat\rho$.

For $\max\{0,k-q\} \leq a \leq \min\{k,p\}$, we define $\hat \tau_a \in \Omega^{n-1}(\overline G)$ to be the coefficient of $\alpha^a \beta^{k-a}$ in the expansion of 
\begin{displaymath}
 \hat \tau_{\alpha,\beta}:=(\alpha \sigma_2+\omega_{2,1}) \wedge \ldots \wedge (\alpha \sigma_{p+1}+\omega_{p+1,1}) \wedge (\beta \sigma_{p+2}+\omega_{p+2,1}) \wedge \ldots \wedge (\beta \sigma_n+\omega_{n,1}).
\end{displaymath}

If $W_2,\ldots,W_n$ are tangent vectors to $\overline G$, then 
$ \hat \tau_{\alpha,\beta}(W_2,\ldots,W_n)$ equals the determinant of the matrix whose entries are $(\alpha \sigma_i+\omega_{i,1})(W_j)$ for $2 \leq i \leq p+1, 2 \leq j \leq n$ and $(\beta \sigma_i+\omega_{i,1})(W_j)$ for $p+2 \leq i \leq n, 2 \leq j \leq n$.
In the exceptional case $2p=2k=n-1$, we define
\begin{equation}\label{omega_especial}
\hat \tau_{ex} := \sigma_{p+2} \wedge \dots \wedge \sigma_n \wedge \omega_{p+2,1} \wedge \dots \wedge \omega_{n,1} \in \Omega^{n-1}(\overline G).
\end{equation}

We denote
\begin{equation}\label{eq_omega_a}
\hat \omega_a := \hat \tau_a \wedge \rho \in \Omega^{n-1+r}(\overline G),\quad \max\{0,k-q\}\leq a\leq\min\{k,p\},
\end{equation}
and, if $n$ is odd, 
\begin{equation}\label{eq_omega_ex}
\hat \omega_{ex}:=\hat \tau_{ex} \wedge \rho \in \Omega^{n-1+r}(\overline G).
\end{equation}

It is obvious that $\hat \omega_a$ and $\hat \omega_{ex}$ are invariant under $H$ and that they vanish on vectors which are tangent to the fiber of $\tilde \Pi$. Hence there are unique forms $\omega_a,\omega_{ex} \in \Omega^{n-1+r}(V \times F(n,p+1))$ with $\tilde \Pi^*\omega_a=\hat \omega_a,\tilde \Pi^*\omega_{ex}=\hat \omega_{ex}$.

From the projection formula \eqref{eq_projection_formula} and the fact that $\tilde\rho$ is the volume form on the fiber of $\tilde{\Pi}$, it follows that for each $\beta \in \Omega^*(V \times F(n,p+1))$
\begin{align}
\hat{\Pi}_{*}(\tilde{\Pi}^{*}\beta \wedge \tilde{\rho}) 
& = \Pi_{*}( \tilde{\Pi}_{*}(\tilde{\Pi}^{*}\beta \wedge \tilde{\rho}) ) \nonumber \\
& = \Pi_{*} (\beta \wedge \tilde{\Pi}_{*}\tilde{\rho}) \nonumber \\
& =\vol(H)  \Pi_{*} \beta, \label{rel_q_qhat}
\end{align}
where we used that $\tilde \Pi_{*}\tilde\rho=\int_H \tilde \rho=\vol H$ (see for instance \cite[(12.11)]{santalo76}).

By~\eqref{eq_fiber_int_gen}, the push-forward of some form $\eta\in\Omega^{n-1+\hat r}(\overline G)$ at the point $(x,v) \in \V \times S^{n-1}$, evaluated at vectors $w_1,\dots,w_{n-1} \in T_{(x,v)}(\V\times S^{n-1})$, is given by 
\begin{equation}\label{eq_fiber_int}
\hat \Pi_*\eta|_{(x,v)}(w_1,\dots,w_{n-1})=\int_{\hat \Pi^{-1}(x,v)}\eta|_{(x,g)}(W_1,\dots,W_{n-1},-),
\end{equation} 
where $g \in G$ is such that $\hat \Pi(x,g)=(x,v)$ and $W_j\in T_{(x,g)}\overline G$ are lifts of $w_j$.

\begin{definition} \label{def_skpi}
Define the smooth flag area measure $S_k^{(p),i}$ by the linear combination
\begin{equation} \label{eq_form_skpi}
\omega:=c_{n,k,p,i} \sum_{a=\min\{k,p\}-m}^{\min\{k,p\}} \binom{\min\{k,p\}-a}{i} \omega_a;
\end{equation}
and, if $n$ is odd and $p=k=\frac{n-1}{2}$, the smooth flag area measure $\tilde S_{\frac{n-1}{2}}^{\left(\frac{n-1}{2}\right)}$  by the form $\omega_{ex}$. 
\end{definition}

\subsection{An integral formula for flag area measures} 

We first recall some notions from \cite{zaehle86}. In the following, we use the convention 
\begin{equation} \label{eq_convention_infty}
\frac{\infty}{\sqrt{1+\infty^2}}=\lim_{\kappa \to \infty} \frac{\kappa}{\sqrt{1+\kappa^2}}=1.
\end{equation}

Let $K$ be a convex body in $\V$, let $x\in\partial K$ and let $(x,v)\in\nc(K)$. 
Then there is a positive orthonormal basis $a_i=a_i(K;x,v),i=2,\ldots,n$ of $v^\perp$ and real numbers $\kappa_i=\kappa_i(K;x,v) \in [0,\infty], i=2,\dots,n$ such that the vectors $w_i:=\frac{1}{\sqrt{1+\kappa_i^2}}(a_i,\kappa_ia_i) \in T_{(x,v)}\nc(K)$ form a positive orthonormal basis of $T_{(x,v)}\nc(K)$. The $\kappa_i$ are called \emph{generalized curvatures}, the $a_i$ are the \emph{generalized curvature directions}. 

The space spanned by the generalized curvature directions with generalized curvature $\kappa \in [0,\infty]$ is unique. If all $\kappa_i$ are finite, we can define a linear operator $S_{x,v}:v^\perp \to v^\perp$ by $S_{x,v}a_i=\kappa_ia_i$. We call it the generalized shape operator. 

We will write formulas involving $S_{x,v}$ even if some of the $\kappa_i$ are infinite. The corresponding term is then to be understood in the sense of a limit. 

If $K$ is smooth and $v=\nu(x)$ is the outer normal vector, then $\kappa_2,\ldots,\kappa_n \in [0,\infty)$ are the principal curvatures of the boundary and the $a_i$ are the principal curvature directions and $S_{x,v}:T_x\partial K\to T_x\partial K$ is the usual shape operator.  

\begin{theorem}  \label{thm_general_body}
Let $k,p,i,\beta$ be as in Theorem \ref{mainthm_construction}. Then for a compact convex body $K$,
\begin{align}
& S_k^{(p),i}(K,\beta)  = c_{n,k,p,i} \binom{n-1}{k} \sum_{a=\min\{k,p\}-m}^{\min\{k,p\}-i} \binom{\min\{k,p\}-a}{i} \binom{k}{a} \times \nonumber \\ 
& \times \int_{\nc(K)} \int_{\Gr_{p+1}(v)} \frac{\mathbf{1}_{(v,E) \in \beta}}{\prod_{i=2}^n\sqrt{1+\kappa_i^2}} D(S_{x,v}[n-k-1],\Pi_E[a],\Pi_{E^\perp}[k-a]) dE d\mathcal{H}^{n-1}(x,v). \label{eq_integral_formula}
\end{align}
In the case $k \leq p$, \eqref{eq_integral_formula} simplifies to 
\begin{align}
& S_k^{(p),i}(K,\beta) = c_{n,k,p,i} \binom{n-1}{k} \binom{k}{i} \times \\ 
& \quad \times \int_{\nc(K)} \int_{\Gr_{p+1}(v)} \frac{\mathbf{1}_{(v,E) \in \beta}}{\prod_{i=2}^n\sqrt{1+\kappa_i^2}} D(S_{x,v}[n-k-1],\Id[k-i],\Pi_{E^\perp}[i]) dE d\mathcal{H}^{n-1}(x,v). \label{eq_integral_formula_simplified}
\end{align}
If $n$ is odd and $p=k=\frac{n-1}{2}$,  
\begin{align*}
 &\tilde S^{\left(\frac{n-1}{2}\right)}_\frac{n-1}{2}(K,\beta) = \\ 
 &(-1)^\frac{n-1}{2} \int_{\nc(K)} \int_{\Gr_{p+1}(v)} \frac{\mathbf{1}_{(v,E) \in \beta}}{\prod_{i=2}^n\sqrt{1+\kappa_i^2}} \det(\Pi_{E^\perp} \circ S_{x,v}:E \cap v^\perp \to E^\perp) dE d\mathcal{H}^{n-1}(x,v).
\end{align*}
\end{theorem}

\proof
Let $K$ be a convex body in $\V$, let $x\in\partial K$ and let $(x,v)\in\nc(K)$. 
We denote by $\kappa_i:=\kappa_i(K;x,v)$, $i=2,\dots,n$, the generalized principal curvatures of $K$ at $(x,v)$ with associated generalized principal directions $a_i:=a_i(K;x,v)$. We order them in such a way that $v,a_2,\ldots,a_n$ is a positive orthonormal basis of $V$.
 
Let $(x,g) \in \overline G$ with $\hat \Pi(x,g)=(x,v)$ and let $g_1,\ldots,g_n$, as before, denote the columns of $g$. Let $E \in \Gr_{p+1}(v)$ be the linear span of $g_1,\ldots,g_{p+1}$. Let $w_i:=\frac{1}{\sqrt{1+\kappa_i^2}}(a_i,\kappa_ia_i) \in T_{(x,v)}\nc(K)$. Then $w_2,\ldots,w_n$ form a positive orthonormal basis of $T_{(x,v)}\nc(K)$. 

By $W_i\in T_{(x,g)}\overline G$ denote a lift of $w_i$, $2\leq i\leq n$, i.e. $d\hat \Pi|_{(x,g)}(W_i)=w_i$.

By definition,  
\begin{equation}\label{eq_ev_sigma}
\sigma_j|_{(x,g)}(W_i)=\langle g_j,d\pi_1(W_i)\rangle=
\frac{1}{\sqrt{1+\kappa_i^2}}\langle g_j,a_i\rangle
\end{equation}
and 
\begin{equation}\label{eq_ev_omega}
\omega_{j,1}|_{(x,g)}(W_i)=
\frac{1}{\sqrt{1+\kappa_i^2}}\langle g_j,\kappa_ia_i\rangle.
\end{equation}

Let $f$ be a smooth function on $F(n,p+1)$. We first compute  
\begin{displaymath}
\Pi_*(f\omega_a)|_{(x,v)}(w_2,\ldots,w_n).
\end{displaymath}

By \eqref{rel_q_qhat}, 
\begin{displaymath}
\Pi_*(f\omega_a)=\frac{1}{\vol H} \hat \Pi_*(\tilde \Pi^* f \wedge \hat \omega_a \wedge \tilde \rho).
\end{displaymath}

Recall that $\rho,\tilde\rho,\hat \rho$ denote the volume forms on the fibers of $\Pi, \tilde\Pi, \hat\Pi$. Using \eqref{eq_fiber_int} and \eqref{eq_omega_a}, we have
\begin{align*}
\Pi_*(f & \omega_a)|_{(x,v)}(w_2,\dots,w_n) \\
&=\frac{1}{\vol H} \hat \Pi_*(\tilde \Pi^{*}f \wedge \hat \omega_a \wedge \tilde\rho)|_{(x,v)}(w_2,\ldots,w_n)\\
&= \frac{1}{\vol H} \int_{\hat \Pi^{-1}(x,v)} \tilde \Pi^*f \cdot (\hat \omega_a \wedge \tilde\rho)|_{(x,g)}(W_2,\dots,W_n,-)\\
&=\frac{1}{\vol H} \int_{\hat \Pi^{-1}(x,v)} \tilde \Pi^{*}f \cdot (\hat \tau_a \wedge \rho \wedge \tilde\rho)|_{(x,g)}(W_2,\dots,W_n,-)\\
&=\frac{1}{\vol H} \int_{\hat \Pi^{-1}(x,v)} \tilde \Pi^{*}f \cdot (\hat \tau_a)|_{(x,g)}(W_2,\ldots,W_n)  \cdot \hat\rho.
\end{align*}
 
Since $\hat \tau_a$ and $\hat\rho$ are invariant under the action of $H$, the expression we are integrating is clearly invariant under the action of $H$ on $\hat \Pi^{-1}(x,v)$. Since we are integrating over the left invariant volume form on $\SO(n-1)$, given as a product of the volume form on $H$ and on $\Gr_{p+1}(v)$ we can interpret the integral as an integral over $\SO(n-1)/H \cong \Gr_{p+1}(v)$. In doing so, we get a factor $\vol(H)$ (see \cite[Theorem 1.48]{sepanski}). Hence 
\begin{align} \label{eq_push_forward_omega_a}
 \Pi_*(f \omega_a)|_{(x,v)}(w_2,\ldots,w_n) & = \int_{\Gr_{p+1}(v)} f(v,E)\cdot \hat \tau_a|_{(x,g)}(W_2,\ldots,W_n) dE.
\end{align}

By definition of $\hat \tau_a$, \eqref{eq_ev_sigma} and \eqref{eq_ev_omega}, $\hat \tau_a|_{(x,g)}(W_2,\ldots,W_n)$ equals the coefficient of $\alpha^a \beta^{k-a}$ in the expansion of the determinant of the matrix $(M_{i,j})_{2 \leq i,j \leq n}$ with
\begin{equation} \label{eq_matrix_m}
M_{i,j}:= \frac{1}{\sqrt{1+\kappa_i^2}}\cdot\begin{cases} 
\alpha\langle g_j,a_i\rangle + \kappa_i \langle g_j,a_i\rangle, & 2 \leq j \leq p+1 \\ 
\beta \langle g_j,a_i\rangle+\kappa_i \langle g_j,a_i\rangle, &  p+2 \leq j \leq n \, .
\end{cases}
\end{equation} 
We define $(n-1) \times (n-1)$-matrices $A, B, C$ by 
\begin{align*}
A_{i,j} & := \begin{cases} 
\langle g_j,a_i\rangle, & 2 \leq j \leq p+1 \\ 
0, &  p+2 \leq j \leq n, \, 
\end{cases}\\
B_{i,j} & := \begin{cases} 
0, & 2 \leq j \leq p+1 \\ 
\langle g_j,a_i\rangle, &  p+2 \leq j \leq n, \, 
\end{cases}\\
C_{i,j} & :=\kappa_i \langle g_j,a_i\rangle, 2 \leq j \leq n.
\end{align*} 
Then 
\begin{displaymath}
M_{i,j}=\frac{1}{\sqrt{1+\kappa_i^2}}(\alpha A+\beta B+C)_{i,j} 
\end{displaymath}
Thus, the coefficient of $\alpha^a\beta^{k-a}$ of the determinant of $M$ equals the mixed discriminant
\begin{displaymath}
\binom{n-1}{k}\binom{k}{a}\frac{1}{\prod_{i=2}^n\sqrt{1+\kappa_i^2}}D(A[a],B[k-a],C[n-k-1]). 
\end{displaymath}

Notice that $A,B,C$ are the matrices of the orthogonal projection to $E\cap v^{\perp}$, the orthogonal projection to $E^{\perp}$, and the generalized shape operator, all computed with respect to the bases $\{a_2,\ldots,a_n\},\{g_2,\ldots,g_n\}$.

Equation \eqref{eq_integral_formula} now follows from Definition \ref{def_skpi}. 

Let us show that \eqref{eq_integral_formula} simplifies to \eqref{eq_integral_formula_simplified} in the case $k \leq p$. Since $\Id=\Pi_E+\Pi_{E^{\perp}}$, we have
\begin{align*}
\binom{k}{i} & D(S_{x,v}[n-k-1], \Id[k-i],\Pi_{E^{\perp}}[i]) \\
& =\sum_{a=0}^{k-i} \binom{k}{i}\binom{k-i}{a}D(S_{x,v}[n-k-1],\Pi_E[a],\Pi_{E^{\perp}}[k-a])\\
& =\sum_{a=0}^{k-i} \binom{k-a}{i}\binom{k}{a} D(S_{x,v}[n-k-1],\Pi_E[a],\Pi_{E^{\perp}}[k-a])\\
 & =\sum_{a=k-m}^{k-i} \binom{k-a}{i}\binom{k}{a}D(S_{x,v}[n-k-1],\Pi_E[a],\Pi_{E^{\perp}}[k-a]),
\end{align*}
where the last line follows from the fact that if $0 \leq a <k-m$, then $k-a>n-p-1$ and the mixed discriminant vanishes since $\Pi_{E^\perp}$ has rank $n-p-1$.

In the exceptional case $p=q=k=\frac{n-1}{2}$, the argument is similar, using 
\begin{displaymath}
 \hat \tau_{ex}(W_2,\ldots,W_n) = \frac{(-1)^{\frac{n-1}{2}}}{\prod_{i=2}^n \sqrt{1+\kappa_i^2}} \det(\Pi_{E^\perp} \circ S_{x,v}:E \cap v^\perp \to E^\perp).
\end{displaymath}
\endproof

\subsection{Proof of Theorems \ref{mainthm_construction} and \ref{mainthm_smooth}} 

\proof[Proof of Theorem \ref{mainthm_smooth}]
The normal cycle of a smooth convex body $K$ is the image of the smooth map $\partial K \to SV, x \mapsto (x,\nu(x))$, where $\nu:\partial K \to S^{n-1}$ is the Gauss map. To transform the integral over the normal cycle into an integral over the boundary, we note that the Jacobian of this map is ${\prod_{i=2}^n\sqrt{1+\kappa_i^2}}$.
\endproof

\proof[Proof of Theorem \ref{mainthm_construction}]
First observe that for $p=0$, $S_k^{(0),0}=S_k$ (the usual surface area measure), which satisfies the formula. Fix $1\leq p, k \leq n-1$. 

Let $x \in \inte F$, where $F$ is a face of $P$ of dimension $\ell$. Let $(x,v) \in \nc(P)$. Every vector tangent to $F$ is a generalized curvature direction with generalized curvature $0$. Every vector in $F^\perp \cap v^\perp$ is  
a generalized curvature direction with generalized curvature $+\infty$. We may therefore choose $a_i:=v_i,i=2,\ldots,n$, where $v=v_1,v_2,\ldots,v_n$ is a positive orthonormal basis of $V$ with $v_2,\ldots,v_{n-\ell}$ spanning $F^\perp \cap v^\perp$ and $v_{n-\ell+1},\ldots,v_n$ spanning $F$. Then $\kappa_2=\ldots=\kappa_{n-\ell}=+\infty$ and $\kappa_{n-\ell+1}=\ldots=\kappa_n=0$. 

Define the matrices $A,B,C$ as in the proof of Theorem \ref{thm_general_body}. Then, 
\begin{displaymath}
\frac{1}{\prod_{i=2}^n\sqrt{1+\kappa_i^2}}D(A[a],B[k-a],C[n-k-1])=D(\tilde A[a],\tilde B[k-a],\tilde C[n-k-1]),
\end{displaymath}
where $\tilde A,\tilde B,\tilde C$ are obtained from $A,B,C$ by multiplying the $i$-th row by $\frac{1}{\sqrt{1+\kappa_i^2}}$. 

The last $\ell$ rows of $\tilde C$ vanish, hence the mixed discriminant vanishes if $\ell>k$.

Similarly, the first $n-\ell-1$ rows in $\tilde A$ and $\tilde B$ vanish, hence the mixed discriminant vanishes if $\ell<k$. 

Let us next consider the case $\ell=k$. The matrix $M$ from \eqref{eq_matrix_m} is then given by
\begin{displaymath}
M_{i,j}:= \begin{cases} 
\langle g_j,v_i\rangle, & 2 \leq i \leq n-k \\ 
\alpha \langle g_j,v_i\rangle, & n-k+1 \leq i \leq n, 2 \leq j \leq p+1  \\ 
\beta \langle g_j,v_i\rangle, & n-k+1 \leq i \leq n, p+2 \leq j \leq n \, .
\end{cases}
\end{displaymath} 
It is easy to see that if we define a matrix $M'$ in an analogous way, but using other orthonormal bases of $E \cap v^\perp, E^\perp, F^\perp \cap v^\perp, F$, then $\det M'=\epsilon \det M$, where $\epsilon=\pm 1$ depends on whether the orientations on $(E \cap v^\perp) \oplus E^\perp$ and $(F^\perp \cap v^\perp) \oplus F$ agree or not.   

We use bases $e_1,\ldots,e_p; g_1,\ldots,g_{n-1-p}; f_1,\ldots,f_k; h_1,\ldots,h_{n-1-k}$ as in Lemma \ref{lem_orbits_grassmann_orthogonal}. Rearranging the rows and columns, the matrix $M'$ has a diagonal block shape, with $m$ blocks of the type $\left(\begin{matrix}
\sin \theta_i' & -\cos \theta_i'\\ \alpha \cos \theta_i' & \beta \sin \theta_i' 
\end{matrix}\right)$,  $\min\{k,p\}-m$ diagonal entries $\alpha$; $\max\{k,p\}-p$ diagonal entries $\beta$, and all other diagonal entries $1$. 

Hence 
\begin{displaymath}
\det M'=\pm \alpha^{\min\{k,p\}-m} \beta^{\max\{k,p\}-p} \prod_{i=1}^m (\alpha \cos^2 \theta_i'+\beta \sin^2 \theta_i') \, .
\end{displaymath}
Since $\det M'=\epsilon$ in the case $\alpha=\beta=1$ (in this case $M'$ is just the transformation matrix between the two bases), we actually have 
\begin{displaymath}
\det M'=\epsilon \alpha^{\min\{k,p\}-m} \beta^{\max\{k,p\}-p} \prod_{i=1}^m (\alpha \cos^2 \theta_i'+\beta \sin^2 \theta_i'),
\end{displaymath}
and hence 
\begin{displaymath}
\det M=\epsilon \det M'= \sigma_{\alpha,\beta}(E,F):= \alpha^{\min\{k,p\}-m} \beta^{\max\{k,p\}-p} \prod_{i=1}^m (\alpha \cos^2 \theta_i'+\beta \sin^2 \theta_i'),
\end{displaymath}
where $\theta_1',\ldots,\theta_m'$ are the principal angles between $E \cap v^\perp$ and $F$.

Let $\sigma^a$ be the coefficient of $\alpha^a \beta^{k-a}$ in $\sigma_{\alpha,\beta}$. Note that $\sigma^a=0$ if $a < \min\{k,p\}-m$ or if $k-a<\max\{k,p\}-p$. We have 
\begin{align*}
\sum_{a=\min\{k,p\}-m}^{\min\{k,p\}} \alpha^a \beta^{k-a} \sigma^a & = \sigma_{\alpha,\beta} \\
& =\alpha^{\min\{k,p\}-m} \beta^{\max\{k,p\}-p} \prod_{i=1}^m (\alpha \cos^2 \theta_i'+\beta \sin^2 \theta_i')\\
& = \alpha^{\min\{k,p\}-m} \beta^{\max\{k,p\}-p} \prod_{i=1}^m (\alpha +(\beta-\alpha) \sin^2 \theta_i').
\end{align*}
Substituting $\alpha:=1, \beta:=t+1$ for some variable $t$, we obtain that  
\begin{displaymath}
\sum_{a=\min\{k,p\}-m}^{\min\{k,p\}} (t+1)^{k-a}\sigma^a=(t+1)^{\max\{k,p\}-p} \prod_{i=1}^m (1+t \sin^2 \theta_i').
\end{displaymath}
We divide both sides by $(t+1)^{\max\{k,p\}-p}$ and obtain
\begin{displaymath}
\sum_{a=\min\{k,p\}-m}^{\min\{k,p\}} (t+1)^{\min\{k,p\}-a}\sigma^a=\prod_{i=1}^m (1+t \sin^2 \theta_i')=\prod_{i=1}^m (1+t \cos^2 \theta_i),
\end{displaymath}
with $\theta_1,\ldots,\theta_m$ being the principal angles between $E^\perp$ and $F$.

Comparing the coefficient of $t^i$ on both sides yields 
\begin{displaymath}
\sum_{a=\min\{k,p\}-m}^{\min\{k,p\}-i} \binom{\min\{k,p\}-a}{i} \sigma^a =\sigma_i(\cos^2 \theta_1,\ldots,\cos^2 \theta_m)=\sigma_i(E^\perp,F).
\end{displaymath}

Taking into account \eqref{eq_push_forward_omega_a} and \eqref{eq_form_skpi} finishes the proof of \eqref{eq_Skpi}. 

Let us finally study the case $2p=2k=n-1$, adapting the argument from above. Define the linear operator $\tilde S_{x,v}$ by 
\begin{displaymath}
\tilde S_{x,v} v_i=\frac{\kappa_i}{\sqrt{1+\kappa_i^2}}v_i=\begin{cases} v_i, & 2 \leq i \leq n-k-1\\0, & n-k \leq i \leq n.\end{cases}
\end{displaymath}
Then, using \eqref{eq_tildesigma_perp} and \eqref{eq_tildesigma_comm}, 
\begin{align*}
\frac{1}{\prod_{i=2}^n\sqrt{1+\kappa_i^2}}  & \det(\Pi_{E^\perp} \circ S_{x,v}:E \cap v^\perp \to E^\perp) = \det(\Pi_{E^\perp} \circ \tilde S_{x,v}:E \cap v^\perp \to E^\perp)\\
& = \det \left(\left\langle \tilde S_{x,v}g_i,g_j\right\rangle\right)_{\substack{2 \leq i \leq p+1 \\p+2 \leq j \leq n}}\\
& = \det \left(\left\langle \tilde S_{x,v} \sum_{a=2}^n \langle g_i,v_a\rangle v_a,\sum_{b=2}^n \langle g_j,v_b\rangle v_b\right\rangle\right)_{\substack{2 \leq i \leq p+1 \\p+2 \leq j \leq n}}\\
& = \det \left(\sum_{a=2}^{n-k-1} \langle g_i,v_a\rangle \cdot  \langle g_j,v_a\rangle\right)_{\substack{2 \leq i \leq p+1 \\p+2 \leq j \leq n}}\\
& = \det B_{F^\perp \cap v^\perp,E \cap v^\perp} \det B_{F^\perp \cap v^\perp,E^\perp}\\
& = \tilde \sigma_{\frac{n-1}{2}}(F^\perp \cap v^\perp,E\cap v^\perp)\\
& = (-1)^{\frac{n-1}{2}}\tilde \sigma_{\frac{n-1}{2}}(E^\perp,F).
\end{align*}

%
\endproof

\subsection{Proof of Theorem \ref{mainthm_properties}} 

\begin{enumerate}
\item[(i)] 
To show Statement (i), it is enough to use the expression for $S_k^{(p)}$ given in Proposition~\ref{prop_explicit_smk} and recall that, as shown in Section \ref{sec_jordan}, we have $\cos(E^{\perp},F)^2=\sigma_{\min\{k,p\}}(E^{\perp},F)$.

\item[(ii)] Translation invariance and homogeneity follow from the corresponding properties of the forms $\omega_a$. 

\item[(iii)] Direct from the fact that the principal angles between a pair of subspaces are invariant under $\O(n)$. 

\item[(iv)] Since the elementary symmetric function of positive numbers is positive, the integrand is positive. 

\item[(v)] 
By a change of variables, Statement (v) follows if, for every $g\in\O(n)$, we have 
\begin{displaymath}
\tilde\sigma_{\frac{n-1}{2}}((gE)^{\perp},gF)=\det g\,\tilde\sigma_{\frac{n-1}{2}}(E^{\perp},F). 
\end{displaymath}

If $g \in \SO(n)$, then $g$ maps positive orthonormal bases to such and the above equation follows. To see what happens for $g \in \O(n) \setminus \SO(n)$, it is enough to look at $g=-\mathrm{Id}$ (which has determinant $-1$ since $n$ is odd). Under this reflection, the orientation of the space $\W=v^{\perp}$ changes and hence the sign of $\tilde\sigma_{\frac{n-1}{2}}(E^{\perp},F)$ changes. 

\item[(vi)]
Let $P$ be a polytope in $\V$ and $\beta \in \mathcal{B}(S^{n-1})$. Then, the right-hand side in \eqref{eq_Skpi} equals 
\begin{displaymath}
c_{n,k,p,i}\sum_{F\in\F_k(P)} \vol_k(F) \int_{\n(P,F)} \left( \int_{\Gr_{p+1}(v)} \sigma_i(E^{\perp},F)dE\right) \mathbf{1}_{v\in\beta} dv. 
\end{displaymath}
The term in brackets is a constant which can be computed with Corollary \ref{cor_expectation_sigma_i}. 

Since $E^\perp$ ranges over an $(n-p-1)$-plane in the $(n-1)$-dimensional space $\W:=v^\perp$, we have $n'=n-1$. Hence, with $m:=\min\{k,n-k-1,p,n-p-1\}$,
\begin{align*}
S_k^{(p),i}(P,\pi^{-1}(\beta)) & = c_{n,k,p,i} \binom{m}{i} \binom{|n-p-k-1|+m}{i}\binom{n-1}{i}^{-1} \\
& \qquad \cdot \sum_{F\in\F_k(P)} \vol_k(F) \int_{\n(P,F)}\mathbf{1}_{v \in \beta}dv\\
& = S_k(P,\beta).
\end{align*}
By approximation, this formula holds for arbitrary convex bodies. 

Finally, to prove that $ \tilde S_\frac{n-1}{2}^{(\frac{n-1}{2})}(K,\pi^{-1}(\beta))=0$, we remark that $(K,\beta) \mapsto  \Phi(K,\beta):=\tilde S_\frac{n-1}{2}^{(\frac{n-1}{2})}(K,\pi^{-1}(\beta))$ satisfies the conditions (B1)-(B5) of Schneider's characterization result \cite[Satz 2]{schneider75b} and is therefore a linear combination of the euclidean area measures. In particular, $\Phi(gK,g\beta)=\Phi(K,\beta)$ for each $g \in \mathrm{O}(n)$. On the other hand, by (iv) we also have $\Phi(gK,g\beta)=\det g\, \Phi(K,\beta)$ for all $g \in \mathrm{O}(n)$. Both equations can hold simultaneously only if $\Phi \equiv 0$. 
\hfill $\square$
\end{enumerate}

\subsection{Proof of Theorem \ref{thm_basis}} 

Since $S_k^{(p),i}$ with $0 \leq i \leq m$ are elements of the space $\FlagArea^{(p),sm,\SO(n)}_k$ and this space is of dimension $m+1$, it remains to prove that these elements are linearly independent. Otherwise, there would be some fixed $k,p$ and constants $c_i$ such that 
\begin{displaymath}
 \sum_{i=0}^m c_i S_k^{(p),i}=0.
\end{displaymath}

Take a polytope of dimension $k$. Let $F$ be its only $k$-face. Fix a unit vector $v$ orthogonal to $F$ and a $(p+1)$-dimensional space $E$ containing $v$. Taking $(f_j)_j$ a sequence of smooth functions on $F(n,p+1)$ with $f_j(v,E)=1$ and whose supports shrink to $(v,E)$, we obtain that 
\begin{displaymath}
0= \sum_{i=0}^m c_i \sigma_i(E^\perp,F)=\sum_{i=0}^m c_i \sigma_i(\cos^2 \theta_1,\ldots,\cos^2 \theta_m),
\end{displaymath}
where $\theta_1,\ldots,\theta_m$ are the Jordan angles between $E^\perp$ and $F$. Since we may choose $E$ arbitrarily, the numbers $\cos^2 \theta_1,\ldots,\cos^2\theta_m$ are arbitrary numbers in $\{0 \leq x_1 \leq \cdots \leq x_m \leq 1\}$. With the $\sigma_i$ being algebraically independent, it follows that each coefficient vanishes. 

This proves that the $S_k^{(p),i}, 0 \leq i \leq m$ are linearly independent. Since these flag area measures are $\mathrm{O}(n)$-covariant, while $\tilde S_\frac{n-1}{2}^{(\frac{n-1}{2})}$ is not $\mathrm{O}(n)$-covariant, there can also be no linear relation involving $\tilde S_\frac{n-1}{2}^{(\frac{n-1}{2})}$.
\hfill $\square$

\section*{Acknowledgements}
The first named author was supported by DFG grants AB 584/1-1 and 584/1-2. The second named author was supported by DFG grant BE 2484/5-2. The authors thank the Oberwolfach Research Institute for Mathematics for its hospitality and support, where this project was started and worked on in various stages of its implementation during several stays at the institute in the framework of the Oberwolfach  Leibniz Fellowship of the third named author. We want to thank the anonymous referee for several suggestions which led to an improvement of the manuscript, in particular Theorem \ref{thm_general_body}.

\def\cprime{$'$}


\begin{thebibliography}{10}

\bibitem{alvarez_fernandes}
Juan Carlos {\'A}lvarez~Paiva and Emmanuel Fernandes.
\newblock Gelfand transforms and {C}rofton formulas.
\newblock {\em Selecta Math. (N.S.)}, 13(3):369--390, 2007.

\bibitem{aomoto87}
Kazuhiko Aomoto.
\newblock Jacobi polynomials associated with {S}elberg integrals.
\newblock {\em SIAM J. Math. Anal.}, 18(2):545--549, 1987.

\bibitem{berline_getzler_vergne}
Nicole Berline, Ezra Getzler, and Mich{\`e}le Vergne.
\newblock {\em Heat kernels and {D}irac operators}.
\newblock Grundlehren Text Editions. Springer-Verlag, Berlin, 2004.
\newblock Corrected reprint of the 1992 original.

\bibitem{bernig_broecker07}
Andreas Bernig and Ludwig Br{\"o}cker.
\newblock {Valuations on manifolds and Rumin cohomology.}
\newblock {\em J. Differ. Geom.}, 75(3):433--457, 2007.

\bibitem{bernig_faifman_opq}
Andreas Bernig and Dmitry Faifman.
\newblock Valuation theory of indefinite orthogonal groups.
\newblock {\em J. Funct. Anal.}, 273(6):2167--2247, 2017.

\bibitem{bernig_fu_hig}
Andreas Bernig and Joseph H.~G. Fu.
\newblock Hermitian integral geometry.
\newblock {\em Ann. of Math.}, 173:907--945, 2011.

\bibitem{bernig_fu_solanes}
Andreas Bernig, Joseph H.~G. Fu, and Gil Solanes.
\newblock Integral geometry of complex space forms.
\newblock {\em Geom. Funct. Anal.}, 24(2):403--492, 2014.

\bibitem{forrester_warnaar}
Peter J.~Forrester and S.~Ole Warnaar.
\newblock The importance of the {S}elberg integral.
\newblock {\em Bull. Amer. Math. Soc. (N.S.)}, 45(4):489--534, 2008.

\bibitem{fu90}
Joseph H.~G. Fu.
\newblock Kinematic formulas in integral geometry.
\newblock {\em Indiana Univ. Math. J.}, 39(4):1115--1154, 1990.

\bibitem{GallagherProulx_1977}
Patrick X.~Gallagher and Ronald J.~Proulx.
\newblock Orthogonal and unitary invariants of families of subspaces.
\newblock In {\em Contributions to algebra (collection of papers dedicated to
  {E}llis {K}olchin)}, pages 157--164. Academic Press, New York, 1977.

\bibitem{goodey_hinderer_hug_rataj_weil}
Paul Goodey, Wolfram Hinderer, Daniel Hug, Jan Rataj, and Wolfgang Weil.
\newblock A flag representation of projection functions.
\newblock {\em Adv. Geom.}, 17(3):303--322, 2017.

\bibitem{hinderer_hug_weil}
Wolfram Hinderer, Daniel Hug, and Wolfgang Weil.
\newblock Extensions of translation invariant valuations on polytopes.
\newblock {\em Mathematika}, 61(1):236--258, 2015.

\bibitem{hug_rataj_weil}
Daniel Hug, Jan Rataj, and Wolfgang Weil.
\newblock Flag representations of mixed volumes and mixed functionals of convex
  bodies.
\newblock {\em J. Math. Anal. Appl.}, 460(2):745--776, 2018.

\bibitem{hug_tuerk_weil}
Daniel Hug, Ines T{\"u}rk, and Wolfgang Weil.
\newblock Flag measures for convex bodies.
\newblock In {\em Asymptotic geometric analysis}, volume~68 of {\em Fields
  Inst. Commun.}, pages 145--187. Springer, New York, 2013.

\bibitem{huybrechts05}
Daniel Huybrechts.
\newblock {\em Complex geometry.}
\newblock Universitext. Springer-Verlag, Berlin, 2005.

\bibitem{james54}
Alan T.~James.
\newblock Normal multivariate analysis and the orthogonal group.
\newblock {\em Ann. Math. Statistics}, 25:40--75, 1954.

\bibitem{jordan1875}
Camille Jordan.
\newblock Essai sur la g\'eom\'etrie \`a {$n$} dimensions.
\newblock {\em Bull. Soc. Math. France}, 3:103--174, 1875.

\bibitem{klain_rota}
Daniel~A. Klain and Gian-Carlo Rota.
\newblock {\em Introduction to geometric probability}.
\newblock Lezioni Lincee. [Lincei Lectures]. Cambridge University Press,
  Cambridge, 1997.

\bibitem{santalo76}
Luis~A. Santal{\'o}.
\newblock {\em Integral geometry and geometric probability}.
\newblock Addison-Wesley Publishing Co., Reading, Mass.-London-Amsterdam, 1976.
\newblock With a foreword by Mark Kac, Encyclopedia of Mathematics and its
  Applications, Vol. 1.

\bibitem{schneider75b}
Rolf Schneider.
\newblock Kinematische {B}er\"uhrma\ss e f\"ur konvexe {K}\"orper.
\newblock {\em Abh. Math. Sem. Univ. Hamburg}, 44:12--23 (1976), 1975.

\bibitem{schneider78}
Rolf Schneider.
\newblock Curvature measures of convex bodies.
\newblock {\em Ann. Mat. Pura Appl. (4)}, 116:101--134, 1978.

\bibitem{schneider_book14}
Rolf Schneider.
\newblock {\em Convex bodies: the {B}runn-{M}inkowski theory}, volume 151 of
  {\em Encyclopedia of Mathematics and its Applications}.
\newblock Cambridge University Press, Cambridge, second expanded edition, 2014.

\bibitem{schneider_weil08}
Rolf Schneider and Wolfgang Weil.
\newblock {\em Stochastic and integral geometry}.
\newblock Probability and its Applications (New York). Springer-Verlag, Berlin,
  2008.

\bibitem{selberg44}
Atle Selberg.
\newblock Remarks on a multiple integral.
\newblock {\em Norsk Mat. Tidsskr.}, 26:71--78, 1944.

\bibitem{sepanski}
Mark R.~Sepanski.
\newblock {\em Compact {L}ie groups}, volume 235 of {\em Graduate Texts in
  Mathematics}.
\newblock Springer, New York, 2007.

\bibitem{spivak2_99}
Michael Spivak.
\newblock {\em A comprehensive introduction to differential geometry. {V}ol.
  {II}}.
\newblock Publish or Perish, Inc., Wilmington, Del., third edition, 1999.

\bibitem{spivak3_99}
Michael Spivak.
\newblock {\em A comprehensive introduction to differential geometry. {V}ol.
  {III}}.
\newblock Publish or Perish, Inc., Wilmington, Del., third edition, 1999.

\bibitem{steiner}
Jakob Steiner.
\newblock {\"U}ber parallele {F}l{\"a}chen.
\newblock {\em Monatsber. Preuß. Akad. Wiss.}, pages 114--118, 1840.
\newblock Ges. Werke, vol. 2, pp. 171--176, Reimer, Berlin, 1882.

\bibitem{wannerer_area_measures}
Thomas Wannerer.
\newblock Integral geometry of unitary area measures.
\newblock {\em Adv. Math.}, 263:1--44, 2014.

\bibitem{zaehle86}
Martina Z{\"a}hle.
\newblock Integral and current representation of {F}ederer's curvature
  measures.
\newblock {\em Arch. Math. (Basel)}, 46(6):557--567, 1986.

\end{thebibliography}

\end{document}